\documentclass[11pt,a4paper]{article}

\usepackage{amsmath,amssymb}
\newcommand{\C}{\mathbb{C}}
\newcommand{\R}{\Omega}

\newcommand{\Z}{\mathbb{Z}}
\def\cA{{\mathcal A}}
\def\cC{{\mathcal C}}
\def\cD{{\mathcal D}}
\def\cS{{\mathcal S}}
\newcommand{\ee}{\varepsilon}
\renewcommand{\div}{{\rm div}\,}
\newcommand{\Supp}{{\rm Supp}\,}
\newcommand{\Int}{\displaystyle \int}
\newcommand{\Frac}{\displaystyle \frac}
\newcommand{\Inf}{\displaystyle \inf}
\newcommand{\Sup}{\displaystyle \sup}
\newcommand{\Lim}{\displaystyle \lim}
\newcommand{\Max}{\displaystyle \max}
\newcommand{\Min}{\displaystyle \min}
\newcommand{\Sum}{\displaystyle \sum}

\def\d{\partial}
\def\ddl{\dot \Delta_l}
\def\ddq{\dot \Delta_q}
\def\tilde{\widetilde}
\def\hat{\widehat}
\newcommand{\D}{\Delta}
\newcommand{\La}{\Lambda}
\newcommand{\n}{\nabla}
\newcommand{\fd}{\frac{d}{2}}
\newcommand{\fdp}{\frac{d}{p}}
\newcommand{\p}{\partial}
\newcommand{\ql}{q_{l}}
\newcommand{\ul}{u_{l}}
\newcommand{\hl}{h_{l}}
\newcommand{\Fl}{F_{l}}
\newcommand{\Gl}{G_{l}}
\newcommand{\Rl}{R_{l}}

\newcommand{\Rlp}{R_{l}'}
\newcommand{\g}{\int_{\Omega^{N}}}

\newcommand{\h}{\hookrightarrow}
\newcommand{\w}{\widetilde{T}}
\newcommand{\de}{\delta}
\newcommand{\del}{\bar{\delta}}
\newcommand{\ka}{\bar{\kappa}}
\newcommand{\cd}{\overline}
\newcommand{\qe}{q_\ee}
\newcommand{\ue}{u_\ee}
\newcommand{\dq}{\delta q}
\newcommand{\du}{\delta u}

\newcommand{\NN}{\frac{N}{p}}
\newcommand{\Om}{\Omega}
\newcommand{\N}{\frac{N}{2}}

\newcommand{\e}{\epsilon}
\newcommand{\va}{\varphi}
\newcommand{\q}{q_{l}}
\newcommand{\ui}{u_{l}}
\newcommand{\T}{\mathbb{T}}

\newtheorem{definition}{Definition}
\newtheorem{theorem}{Theorem}

\newtheorem{proposition}{Proposition}
\newtheorem{corollaire}{Corollary}
\newtheorem{remarka}{Remark}

\newtheorem{lemme}{Lemma}

\def\thesection{\arabic{section}}
\def\theequation{\arabic{section}.\arabic{equation}}
\def\thetheorem{\arabic{section}.\arabic{theorem}}
\def\theproposition{\arabic{section}.\arabic{proposition}}
\def\thecorollary{\arabic{section}.\arabic{corollaire}}
\def\thedefinition{\arabic{section}.\arabic{definition}}
\def\theremark{\arabic{section}.\arabic{remarque}}
\newcommand{\reset}{\setcounter{equation}{equation}\setcounter{theorem}{theorem}
\setcounter{proposition}{0}\setcounter{corollary}{corollaire}
\setcounter{remark}{remark}}

\title{New entropy for Korteweg's system, existence of global weak solution and Prodi-Serrin theorem}
\author{Boris Haspot \thanks{Basque Center of Applied Mathematics, Bizkaia Technology Park, Building 500,
E-48160, Derio (Spain), haspot@bcamath.org }}
\date{}
\begin{document}
\maketitle
\begin{abstract}
This work is devoted to prove new entropy estimates for a general
isothermal model of capillary fluids derived by J.E Dunn and
J.Serrin (1985)  (see \cite{fDS}), which can be used as a phase transition model. More precisely we will derive new estimates for the density and we will give a new structure for the Korteweg system which allow us to obtain the existence of global weak solution. The key of the proof comes from the introduction of a new effective velocity.The proof is widely inspired from the works of A. Mellet and A. Vasseur  (see \cite{fMV}). In a second part, we shall give a  Prody-Serrin blow-up criterion for this system which  widely improves the results of \cite{Hprepa} and the known results on compressible systems.
\end{abstract}
\section{Introduction}
We are concerned with compressible fluids endowed with internal
capillarity. The model we consider  originates from the XIXth
century work by Van der Waals and Korteweg \cite{VW,fK} and was
actually derived in its modern form in the 1980s using the second
gradient theory, see for instance \cite{fDS,fJL,fTN}. The first investigations begin with the Young-Laplace theory which claims that the phases are separated by a hypersurface and that the jump in the pressure across the hypersurface is proportional to the curvature of the hypersurface. The main difficulty consists in describing the location and the movement of the interfaces.\\
Another major problem is to understand whether the interface behaves as a discontinuity in the state space (sharp interface) or whether the phase boundary corresponds to a more regular transition (diffuse interface, DI).
The diffuse interface models have the advantage to consider only one set of equations in a single spatial domain (the density takes into account the different phases) which considerably simplifies the mathematical and numerical study (indeed in the case of sharp interfaces, we have to treat a problem with free boundary).\\
Another approach corresponds to determine equilibrium solutions which classically consists in the minimization of the free energy functional.
Unfortunately this minimization problem has an infinity of solutions, and many of them are physically wrong (some details are given later). In order to overcome this difficulty, Van der Waals in the XIX-th century was the first to add a term of capillarity to select the physically correct solutions, modulo the introduction of a diffuse interface. This theory is widely accepted as a thermodynamically consistent model for equilibria. Alternatively, another way to penalize the high density variations consists in applying a zero order but non-local operator to the density gradient \cite{9Ro}, \cite{5Ro}, \cite{Rohdehdr}.  We refer for a mathematical analysis on this system to \cite{CH,Has1,Has5,Has4}.
Korteweg-type models are based on an extended version of
nonequilibrium thermodynamics, which assumes that the energy of the
fluid not only depends on standard variables but also on the
gradient of
the density.\\
Let us now consider a fluid of density $\rho\geq 0$, velocity field $u\in\R^{N}$ (both are defined on a subset $\Om$ with $\Om=\R^{N}$ or the torus $\T^{N}$), we are now interested in the following
compressible capillary fluid model, which can be derived from a Cahn-Hilliard like free energy (see the
pioneering work by J.- E. Dunn and J. Serrin in \cite{fDS} and also in
\cite{fA,fC,fGP,HM}).
The conservation of mass and of momentum write:
\begin{equation}
\begin{cases}
\begin{aligned}
&\frac{\p}{\p t}\rho+{\rm div}(\rho u)=0,\\
&\frac{\p}{\p t}(\rho u)+{\rm div}(\rho
u\otimes u)-\rm div(\mu\rho\,\n u)-\rm div(\alpha \rho\,\n u^{t})+\n (a\rho^{\gamma})={\rm div}K,
\end{aligned}
\end{cases}
\label{3systeme}
\end{equation}
where the Korteweg tensor read as following:
\begin{equation}
{\rm div}K
=\n\big(\rho\kappa(\rho)\D\rho+\frac{1}{2}(\kappa(\rho)+\rho\kappa^{'}(\rho))|\n\rho|^{2}\big)
-{\rm div}\big(\kappa(\rho)\n\rho\otimes\n\rho\big).
\label{divK}
\end{equation}
$\kappa$ is the coefficient of capillarity and is a regular function of the form $\kappa(\rho)=\kappa\rho^{\alpha}$ with $\alpha\in\R$. The term
${\rm div}K$  allows to describe the variation of density at the interfaces between two phases, generally a mixture liquid-vapor. $P=a\rho^{\gamma}$ with $\gamma\geq 1$ is a general barotropic pressure term.
$\mu$ and $\alpha$ are the two Lam\'e viscosity coefficients and satisfying:
$$\mu>0\;\;\mbox{and}\;\;\mu>\alpha\geq 0.$$
In particular, it allows to write the diffusion tensor on the form $(\mu-\alpha){\rm div}(\rho\n u)+\alpha{\rm div}(\rho Du)$ with $Du =\n u+^{t}\n u$ the strain tensor which implies some energy inequality. More precisely when we multiply the momentum equation by $u$, we have:
\begin{equation}
\begin{aligned}
&\int_{\Om}\big(\rho(t,x)|u(t,x)|^{2}+\frac{a}{\gamma-1}\rho^{\gamma}(t,x)+\kappa|\n\sqrt{\rho}|^{2}(t,x)\big)dx\\
&\hspace{2cm}+\int^{t}_{0}\int_{\Om}\big((\mu-\alpha)\rho(t,x)|\n u|^{2}(t,x)+\alpha \rho(t,x)|D u|^{2}(t,x)\\
&\hspace{3cm}\leq C(\int_{\Om}\big(\rho_{0}(x)|v_{0}(x)|^{2}+\frac{1}{\gamma-1}\rho_{0}^{\gamma}(x)+\kappa|\n\sqrt{\rho_{0}}|^{2}\big)dx).
\end{aligned}
\label{ineg1}
\end{equation}
We now want to recall the existing results on the  existence of global weak solutions on classical compressible Navier Stokes equations and on Korteweg system.
\subsubsection*{Weak solutions}
We wish to prove existence and uniqueness results for $(NHV)$ in
functions spaces very close to energy spaces. In the non isothermal
non capillary case and $P(\rho)=a\rho^{\gamma}$, with $a>0$ and
$\gamma>1$, P-L. Lions in \cite{fL2} proved the global existence of
variational solutions $(\rho,u,\theta)$ to $(NHV)$ with $\kappa=0$
for $\gamma> \frac{N}{2}$ if $N\geq 4$, $\gamma\geq \frac{3N}{N+2}$
if $N=2,3$ and initial data $(\rho_{0},m_{0})$ such that:
$$\Pi(\rho_{0})-\Pi(\bar{\rho}),\;\;\frac{|m_{0}|^{2}}{\rho_{0}}\in
L^{1}(\R^{N}),\;\;\mbox{and}\;\;\rho_{0}\theta_{0}\in
L^{1}(\R^{N}).$$ These solutions are weak solutions in the classical
sense for the equation of mass conservation and for
the equation of the momentum. Notice that the main difficulty for proving Lions' theorem consists
in exhibiting strong compactness properties of the density $\rho$ in
$L^{p}_{loc}$ spaces required to pass to the limit in the pressure
term $P(\rho)=a\rho^{\gamma}$.\\
Let us mention that Feireisl in \cite{fF}  generalized the result to
$\gamma>\frac{N}{2}$ in establishing that we can obtain renormalized
solution without imposing that $\rho\in L^{2}_{loc}$, for this he
introduces the concept of oscillation defect measure evaluating the
lost of compactness.\\
We can finally cite the result from Bresch-Desjardins
in \cite{5BD},\cite{5BD1} where they show the existence of global weak
solution for $(NHV)$ with $\kappa=0$ and with a cold pressure.  In fact they are working with specific viscosity coefficients which verify an algebraic relation. It allows them to
get good estimate on the density by using new energy inequality and to
treat by compactness all the delicate terms as the pressure. In \cite{fMV}, Mellet and Vasseur improve the results of Bresch,Desjardins by dealing with the case of a general pressure $P(\rho)=a\rho^{\gamma}$ with $\gamma>1$.\\
\\
In the case $\kappa>0$, we remark then that the density belongs to
$L^{\infty}(0,\infty,{\dot {H}}^{1}(\R^{N}))$. Hence, in contrast to
the non capillary case one can easily pass to the limit in the
pressure term. However let us emphasize at this point that the above
a priori bounds do not provide any $L^{\infty}$ control on the
density from below or from above. Indeed, even in dimension $N=2$,
$H^{1}$ functions are not necessarily locally bounded. Thus, vacuum
patches are likely to form in the fluid in spite of the presence of capillary forces, which are expected to smooth out the density. It explains why it is so difficult to obtaining the existence of global strong solution in dimension $N=2$. Indeed in order to prove the existence of global weak solution the main difficulty consists in dealing with the quadratic terms in gradient of the density appearing in the capillary tensor. Recently D. Bresch, B. Desjardins and C-K. Lin in \cite{5BDL} got
some global weak solutions for the isotherm Korteweg model with some
specific viscosity coefficients. Indeed, they assume that
$\mu(\rho)=C\rho$ with $C>0$ and $\lambda(\rho)=0$. By choosing
these specific coefficients they obtain a gain of derivatives on the
density $\rho$ where $\rho$ belongs to $L^{2}(H^{2})$.
It is easy at that time with this kind of estimate on the density $\rho$ to get strong convergence on the term of capillarity. However a new difficulty appears
concerning the loss of information on the gradient of $u$ when vacuum existing (indeed the viscosity coefficients are degenerated). It becomes
involved to pass to the limit in the term $\rho_{n} u_{n}\otimes u_{n}$. That's why
the solutions of D. Bresch, B. Desjardins and C-K. Lin require some
specific test
functions which depend on the density $\rho$. \\
In \cite{fH2}, we improve the results of R. Danchin and
B. Desjardins in \cite{fDD} and D. Bresch, B. Desjardins an C-K. Lin in \cite{5BDL} by showing the existence of global weak solution with small initial data in the energy space
for specific choices of the capillary coefficients and with general viscosity coefficient.
Comparing with the results of \cite{5BDL}, we get global weak solutions with general test function $\va\in C^{0}_{\infty}(\R^{N})$not depending on the density $\rho$. In fact we have extracted of the structure of capillarity term a new energy inequality using fractionary derivative which allows a
gain of derivative on the density $\rho$.
\subsection{Derivation of the models}
We are going to prove that we can derive new entropy estimates when we choose specific coefficients for the viscosity and the capillarity. In the sequel we will
consider the following physical coefficients:
$$\mu(\rho)=\mu\rho\;\;\;\mbox{and}\;\;\;\kappa(\rho)=\frac{\kappa}{\rho},$$
with $\mu,\kappa>0$. By computation, we show that in this case, we obtain (see the appendix for more details):
$${\rm div}K=\kappa{\rm div}(\rho\n\n\ln\rho)=\kappa{\rm div}(\rho\n(D\ln\rho)).$$
It is now clear that we have to work with the new variable $v=u+\frac{\kappa}{\mu}\n\ln\rho$. We now want to rewrite system (\ref{3systeme}) in terms of the variables $(\ln\rho,v)$. We have then by considering the mass equation:
$$\rho\p_{t}\ln\rho+\rho u\cdot\n \n\ln\rho+\rho\n\ln\rho\cdot\n u^{t}+\rho\n{\rm div}u=0.$$
Then we obtain the following new system when $\alpha=\frac{\kappa}{\mu}$:
\begin{equation}
\begin{cases}
\begin{aligned}
&\p_{t}\rho+{\rm div}(\rho v)-\frac{\kappa}{\mu}\D\rho=0,\\
&\rho\p_{t}v +\rho u\cdot\n v-\rm div(\mu\rho\,\n v)+\n P(\rho)=0,
\end{aligned}
\end{cases}
\label{3systeme1}
\end{equation}
When $\alpha=0$ and $\kappa=\mu^{2}$, we obtain the following simplified model:
\begin{equation}
\begin{cases}
\begin{aligned}
&\p_{t}\rho+{\rm div}(\rho v)-\frac{\kappa}{\mu}\D\rho=0,\\
&\rho\p_{t}v +\rho u\cdot\n v-\rm div(\mu\rho\,\n v)+\n P(\rho)=0,
\end{aligned}
\end{cases}
\label{3systeme2}
\end{equation}
For more details on the computation, we refer to the appendix. Our goal is now to prove new entropy inequalities for these two systems and to obtain the existence of global weak solution for these two systems. We also prove the existence of global weak solution for the shallow-water system, indeed A. Mellet and A. Vasseur in \cite{fMV} have obtained the stability of global weak solutions for shallow-water system but to get the existence of global weak solution, it remains to construct approximate solutions of the shallow water system conserving the entropy inequalities used in \cite{fMV}. This is exactly what we will show by studying the system (\ref{3systeme1}).
\section{Notations and main result}
\subsection{Existence of global weak solution for Korteweg system}
We say that $(\rho,v)$ is a weak solution of (\ref{3systeme1}) on $[0,T]\times\Omega$, which the following initial conditions
\begin{equation}
\rho_{/t=0}=\rho_{0}\geq0,\;\;\rho u_{/t=0}=m_{0}.
\label{7}
\end{equation}
with:
\begin{equation}
\begin{aligned}
&\rho_{0}\in L^{\gamma}(\Om)\cap L^{1}(\Om),\;\sqrt{\rho_{0}}\n\ln\rho_{0}\in L^{2}(\Om),\;\rho_{0}\geq 0,\\
&\sqrt{\rho_{0}}v_{0}\in L^{2}(\Om),\;\rho_{0}^{\frac{1}{2+\delta}}v_{0}\in L^{2+\delta}(\Om)\;\;\;\mbox{for somme small $\delta$}.
\end{aligned}
\label{8}
\end{equation}
if
\begin{itemize}
\item $\rho\in L^{\infty}_{T}(L^{1}(\Om)\cap L^{\gamma}(\Om))$, $\sqrt{\rho}\in L_{T}^{\infty}(L^{2}(\Om))$,
\item $\sqrt{\rho}v\in L^{\infty}_{T}(L^{2}(\Om))$,
\item $\sqrt{\rho}\, \n v\in L^{2}((0,T)\times\Om)$,
\end{itemize}
with $\rho\geq 0$ and $(\rho,\sqrt{\rho}v)$ satisfying in ${\cal D}^{'}([0,T]\times\R^{N})$:
$$
\begin{cases}
\begin{aligned}
&\p_{t}\rho+{\rm div}(\sqrt{\rho}\sqrt{\rho}v)-\frac{\kappa}{\mu}\D\rho=0,\\
&\rho(0,x)=\rho_{0}(x).
\end{aligned}
\end{cases}
$$
and if the following equality holds for all $\va(t,x)$ smooth test function with compact support such that $\va(T,\cdot)=0$:
\begin{equation}
\begin{aligned}
&\int_{\Om}(\rho v)_{0}\cdot\va(0,\cdot)dx+\int^{T}_{0}\int_{\Om}\sqrt{\rho}(\sqrt{\rho}v)\p_{t}\va+
\sqrt{\rho}u\otimes\sqrt{\rho}v:\n\va dx\\
&\hspace{4cm}+\int^{T}_{0}\int_{\Om}\rho^{\gamma}{\rm div}\va-<\mu\rho\, \n v,\n \va>=0,
\end{aligned}
\label{equadistrib}
\end{equation}
where we give sense to the diffusion terms by rewriting him according to $\sqrt{\rho}$ and $\sqrt{\rho}v$:
$$
\begin{aligned}
<\rho\, \n v,\n\va>=&-\int\sqrt{\rho}(\sqrt{\rho}v_{j})\p_{ii}\va_{j}dx\,dt-\int
2\sqrt{\rho}v_{j}\p_{i}\sqrt{\rho}\p_{i}\va_{j}dx\,dt\\
\end{aligned}
$$
Similarly we have exactly the same type of definition for system (\ref{3systeme2}).
$$
\begin{aligned}
&<\rho\n v,\n\va>=-\int\sqrt{\rho}(\sqrt{\rho}v_{j})\p_{ii}\va_{j}dx\,dt-\int
2\sqrt{\rho}v_{j}\p_{i}\sqrt{\rho}\p_{i}\va_{j}dx\,dt.
\end{aligned}
$$
\subsubsection*{Main results}
We obtain in this paper the existence of global weak solutions (more exactly the stability of global weak solutions) for systems (\ref{3systeme1}) and
(\ref{3systeme2}). For system (\ref{3systeme1}) we obtain the following first theorem.
\begin{theorem}
\label{2.1}
Assume that $\gamma>1$. Let $(\rho_{n},v_{n})$ with $v_{n}=u_{n}+\frac{\kappa}{\mu}\n\log\rho_{n}$ be a sequence of weak solutions of system \ref{3systeme1} satisfying entropy inequalities (\ref{21}) and (\ref{22}), with initial data:
$$(\rho_{n})_{/t=0}=\rho_{0}^{n}(x)\;\;\;\mbox{and}\;\;\;(\rho_{n}v_{n})_{/t=0}=\rho_{0}^{n}v^{n}_{0}(x)$$
where $\rho_{0}^{n}$ and $v^{n}_{0}$ such that:
\begin{equation}
\rho_{0}^{n}\geq0,\;\;\rho_{0}^{n}\rightarrow\rho_{0}\;\;\mbox{in}\;L^{1}(\Omega),\;\;\rho_{0}^{n}v^{n}_{0}\rightarrow\rho_{0}v_{0}\;\;\mbox{in}\;L^{1}(\Omega),
\label{15}
\end{equation}
and satisfy the following bounds (with $C$ constant independent on $n$):
\begin{equation}
\int_{\Om}\big(\rho_{0}^{n}\frac{|v_{0}^{n}|^{2}}{2}+\frac{a(\rho_{0}^{n})^{\gamma}}{\gamma-1}\big)<C,\;\;
\int_{\Om}\frac{1}{\rho_{0}^{n}}|\n \ln \rho_{0}^{n}|^{2}dx<C,
\label{16}
\end{equation}
and:
\begin{equation}
\int_{\Om}\rho_{0}^{n}\frac{|v_{0}^{n}|^{2+\alpha}}{2}dx<C,
\label{17}
\end{equation}
Then, up to a subsequence, $(\rho_{n},\sqrt{\rho_{n}}v_{n},\sqrt{\rho_{n}}u_{n})$ converges strongly to a weak solution $(\rho,\sqrt{\rho}v,\sqrt{\rho}u)$ of (\ref{3systeme1})  satisfying entropy inequalities (\ref{21}) and (\ref{22}) (the density $\rho_{n}$ converges strongly in $C^{0}((0,T), L^{\frac{3}{2}}_{loc}(\Om))$, $\sqrt{\rho_{n}}v_{n}$ converges strongly in $L^{2}(0,T,L^{2}_{loc})$ and the momentum $m_{n}=\rho_{n}v_{n}$ converges strongly in $L^{1}(0,T,L^{1}_{loc}(\Om))$, for any $T>0$).
\end{theorem}
\begin{remarka}
The proof will be strongly inspired of the works of A. Mellet and A. Vasseur in \cite{fMV}.
\end{remarka}
We now obtain exactly thesame result for system (\ref{3systeme2}).
\begin{theorem}
\label{2.2}
Assume that we have a sequence $(\rho_{n},v_{n})$  with $v_{n}=u_{n}+\frac{\kappa}{\mu}\n\log\rho_{n}$ of weak solutions of  system (\ref{3systeme2}) satisfying entropy inequalities (\ref{21a}) and (\ref{22a}) with initial data:
$$(\rho_{n})_{/t=0}=\rho_{0}^{n}(x)\;\;\;\mbox{and}\;\;\;(\rho_{n}v_{n})_{/t=0}=\rho_{0}^{n}v^{n}_{0}(x)$$
where $\rho_{0}^{n}$ and $v^{n}_{0}$ such that:
\begin{equation}
\rho_{0}^{n}\geq0,\;\;\rho_{0}^{n}\rightarrow\rho_{0}\;\;\mbox{in}\;L^{1}(\Omega),\;\;\rho_{0}^{n}v^{n}_{0}\rightarrow\rho_{0}v_{0}\;\;\mbox{in}\;L^{1}(\Omega),
\label{15a}
\end{equation}
and satisfy the following bounds (with $C$ constant independent on $n$):
\begin{equation}
\int_{\Om}\big(\rho_{0}^{n}\frac{|v_{0}^{n}|^{2}}{2}+\frac{a(\rho_{0}^{n})^{\gamma}}{\gamma-1})<C,\;\;
\int_{\Om}\frac{1}{\rho_{0}^{n}}|\n \ln \rho_{0}^{n}|^{2}dx<C,
\label{16a}
\end{equation}
and:
\begin{equation}
\int_{\Om}\rho_{0}^{n}\frac{|v_{0}^{n}|^{2+\delta}}{2}dx<C,
\label{17a}
\end{equation}
Then, up to a subsequence, $(\rho_{n},\sqrt{\rho_{n}}v_{n},\sqrt{\rho_{n}}u_{n})$ converges strongly to a weak solution $(\rho,\sqrt{\rho}v,\sqrt{\rho}u)$ of (\ref{3systeme2}) satisfying entropy inequalities (\ref{21a}) and (\ref{22a}) (the density $\rho_{n}$ converges strongly in $C^{0}((0,T), L^{\frac{3}{2}}_{loc}(\Om))$, $\sqrt{\rho_{n}}v_{n}$ converges strongly in $L^{2}(0,T,L^{2}_{loc})$ and the momentum $m_{n}=\rho_{n}v_{n}$ converges strongly in $L^{1}(0,T,L^{1}_{loc}(\Om))$, for any $T>0$).
\end{theorem}
In the specific case of the system (\ref{3systeme2})  we obtain new blow-up criterion in the case of the torus $\T^{N}$ which improves the results in \cite{Hprepa}. In the sequel, we will set $m=\rho v$ and $q{'}=\rho-1$. We have then the following theorem:
\begin{theorem}
Let  $P$ be a suitably smooth function of the density and  $1\leq p <+\infty$.
Let $m_{0}\in \dot{B}^{\frac{N}{p}-1+\e^{'}}_{p,\infty}$ with $\e^{'}>0$ and $q^{'}_{0}\in \dot{B}^{\NN+\e^{'}}_{p,\infty}$ such that $\rho_{0}\geq c>0$. \\
There exists then a positive time $T$ such that system (\ref{3systeme2}) has a unique solution
$(q^{'},m)$ with $\rho$ bounded away from $0$ and:
$$q^{'}\in \widetilde{C}([0,T],\dot{B}^{\NN+\e^{'}}_{p,1}
)\cap \widetilde{L}^{1}_{T}(\dot{B}^{\NN+2+\e^{'}}_{p,1}),\;\;m\in \widetilde{C}([0,T];\dot{B}^{\frac{N}{p}-1+\e^{'}}_{p,1})\cap\widetilde{L}^{1}([0,T],
\dot{B}^{\NN+1+\e^{'}}_{p,1}).$$
We assume now that $P(\rho)=a\rho$ with $a>0$. Furthermore if $v_{0}\in L^{\infty}$, $\rho_{0}\in B^{1}_{p,\infty}$ for any $1\leq p<+\infty$, $\frac{1}{\rho_{0}}1_{\{|\rho_{0}|\leq\delta\}}\in L^{\infty}(\T^{N})\cap L^{1}(\T^{N})$  (with $0<\delta<1$) and the initial data are in the energy space, it means:
$$\sqrt{\rho_{0}}u_{0}\in L^{2},\,\n\sqrt{\rho_{0}}\in L^{2}\;\;\;\mbox{and}\;\;\;\Pi(\rho_{0})\in L^{1}.$$
(we refer to the proposition \ref{BD} for the definition of $\Pi$) we can then extend the solution beyond $(0,T)$ if:
\begin{equation}
v\in L^{p}_{T}(L^{q}(\T^{N}))\;\;\;\;\mbox{with}\;\;\;\frac{1}{p}+\frac{N}{2q}=\frac{1}{2}\;\;\mbox{and}\;\;1\leq p<+\infty,
\end{equation}
or if for any $\e>0$ arbitrary small:
\begin{equation}
\frac{1}{\rho^{\e}}1_{\{|\rho|\leq\delta\}} \in L^{\infty}_{T}(L^{1}(\T^{N})).\label{blow1}
\end{equation}
\label{theo4}
\end{theorem}
\begin{remarka}
This result is really to consider as a Prodi-Serrin theorem on the effective velocity $v$. In terms of blow-up condition, he improves widely \cite{Hprepa}.In fact the second condition could be improved as follows (we refer to the proof of theorem \ref{theo4}):
\begin{equation}
\frac{1}{\rho^{\e}}1_{\{|\rho|\leq\delta\}} \in L^{p}_{T}(L^{q}(\T^{N}))\;\;\;\;\mbox{with}\;\;\;\frac{1}{p}+\frac{N}{2q}=\frac{1}{2}.\label{blow1}
\end{equation}
\end{remarka}
\begin{remarka}
We could probably extend this previous result for more general pressure terms. However it would requires additional informations on the integrability of the density or of the vacuum, i.e $\frac{1}{\rho}1_{\{|\rho|\leq\delta\}} $.\\
It would be also possible to deal with the euclidian space $\mathbb{\R}^{N}$ (it does not change a lot, except that sometime, we shall use interpolation so that we would need of additional information in terms of low frequencies.
\end{remarka}
\section{New entropies}
\subsubsection{Entropy for the system (\ref{3systeme1}) }
We now want to establish new entropy inequality for system (\ref{3systeme1}) and (\ref{3systeme2}).
We now obtain  the following proposition when $(\rho,u)$ are exact solutions of system (\ref{3systeme1}).
\begin{proposition}
Assume that $(\rho,u)$ are exact solutions of system (\ref{3systeme1}) with $P(\rho)=a\rho^{\gamma}$ ($\gamma\geq 1$) then for all $t>0$:
\begin{equation}
\begin{aligned}
&\int_{\Om}\big[\rho|u|^{2}(t,x)+\kappa|\n\sqrt{\rho}|^{2}(t,x)+\Pi(\rho)(t,x)\big]\,dx+\int_{0}^{t}\int_{\Om}\n\ln\rho\cdot\n\rho^{\gamma}\,dx dt\\
&+(\mu-\gamma)\int^{t}_{0}\int_{\Om}\rho|\n u|^{2}dxdt+\gamma\int^{t}_{0}\int_{\Om}\rho|D u|^{2}dxdt+\kappa\int^{t}_{0}\int_{\Om}\rho(\p_{ij}\ln\rho)^{2}(t,x)dtdx \\
&\hspace{4cm}\leq C( \int_{\R^{N}}\big(\rho_{0}|v_{0}|^{2}(x)+\Pi(\rho_{0}(x))+\kappa|\n\sqrt{\rho}_{0}|^{2}(t,x\big)\,dx).
\end{aligned}
\label{21}
\end{equation}
with $\Pi(s)=s\int^{s}\int^{s}_{0}\frac{P(z)}{z^{2}}dz$.
\label{BD}
\end{proposition}
{\bf Proof:}
We now want to obtain this new entropy by two different ways, one which is very direct and an other using the BD entropy (see \cite{5BD2}).
\subsubsection*{Direct energy inequality via system (\ref{3systeme1})}
When we multiply the momentum equation in (\ref{3systeme1}) by $v$, we obtain:
$$
\begin{aligned}
&\int_{\Om}\big(\rho(t,x)|v(t,x)|^{2}+\Pi(\rho)(t,x)\big)dx+\int^{t}_{0}\int_{\Om}\big(\mu\rho(t,x)|\n v|^{2}(t,x)\\
&\hspace{1cm}+\frac{\kappa}{\mu}P^{''}(\rho)|\n\rho|^{2}(t,x)\big)dtdx\leq C(\int_{\Om}\big(\rho_{0}(x)|v_{0}(x)|^{2}+\Pi(\rho_{0})(x)\big)dx).
\end{aligned}
$$
By the previous inequality and (\ref{ineg1}) we obtain the desired result.
\subsubsection*{BD entropy via system (\ref{3systeme})}
By an other way we obtain the new BD entropy of proposition (\ref{BD})  by multiplying the momentum equation by $\n\ln\rho$. Indeed the only difference with the classical BD entropy (see \cite{5BD2}) is the extra capillarity terms and the term $\frac{\kappa}{\mu}{\rm div}(\rho\n u^{t})$ in (\ref{3systeme}). We have only to treat these two extra terms.
$$
\begin{aligned}
&\int^{t}_{0}\int_{\Om}{\rm div}K\cdot\n\ln\rho\,dxdt=\kappa\int^{t}_{0}\int_{\Om}{\rm div}(\rho\n\n\ln\rho)\cdot\n\ln\rho\,dxdt,\\
&\hspace{3cm}=-\kappa \int^{t}_{0}\int_{\Om}\rho(\p_{i,j}\ln\rho)^{2} dtdx.
\end{aligned}
$$
$$
\begin{aligned}
\int^{t}_{0}\int_{\Om}\frac{\kappa}{\mu}{\rm div}(\rho\n u^{t})\cdot\n\ln\rho\, dtdx&=-\int^{t}_{0}\int_{\Om}\frac{\kappa}{\mu}\rho\p_{j}u_{i}\p_{ij}\ln\rho\, dtdx,\\
&\leq\frac{\kappa}{\mu}(\frac{1}{2\e}\int^{t}_{0}\int_{\Om}\rho (\p_{ij}\ln\rho)^{2}\, dtdx+\frac{\e}{2}\int^{t}_{0}\int_{\Om}\rho|\n u^{t}|^{2}dx.
\end{aligned}
$$
By bootstrap we conclude. For more details on the BD entropy in our case, we refer to the appendix.
{\hfill $\Box$}\\
\\
The following proposition comes from \cite{fMV1}.
\begin{proposition}
Smooth solutions of system (\ref{3systeme1}) satisfy the following inequality when $P(\rho)=a\rho^{\gamma}$ with $\gamma\geq 1$:
\begin{equation}
\begin{aligned}
&\frac{d}{dt}\int_{\Om}\rho\frac{|v|^{2+\de}}{2+\de}+\frac{\nu}{4}\int_{\Om} \rho|v|^{\de}|\n v|^{2}dx\leq
\biggl(\int_{\Om}\big(\rho^{2\gamma-1-\frac{\delta}{2}}\big)^{\frac{2}{2-\de}}dx\biggl)^{\frac{2}{2-\de}}
\big(\int_{\Om}\rho|v|^{2}dx\big)^{\frac{\de}{2}},
\end{aligned}
\label{22}
\end{equation}
for $\delta\in(0,2)$.
\end{proposition}
{\bf Proof:} The proof follows exactly the same lines than in \cite{fMV}.
\subsubsection{Entropy for the system (\ref{3systeme2}) }
By proceeding similarly we obtain the following propositions for system (\ref{3systeme2}).
\begin{proposition}
Assume that $(\rho,u)$ are exact solutions of system (\ref{3systeme2}) then for all $t>0$:
\begin{equation}
\begin{aligned}
&\int_{\Om}\big(\rho|v|^{2}(t,x)+\Pi(\rho(t,x))\big)\,dx+\int_{0}^{t}\int_{\Om}\n\ln\rho\cdot\n\rho^{\gamma}\,dx dt\\
&\hspace{1cm}+\kappa\int^{t}_{0}\int_{\Om}\rho(\p_{ij}\ln\rho)^{2}(t,x)dtdx \leq C( \int_{\Om}\big(\rho_{0}|v_{0}|^{2}(x)+\Pi(\rho_{0}(x))\big)\,dx).
\end{aligned}
\label{21a}
\end{equation}
\label{BD1}
\end{proposition}
{\bf Proof:} The proof follows the same lines as in proposition \ref{BD}.
{\hfill $\Box$}\\
\\
For the following proposition, we refer to \cite{fMV}.
\begin{proposition}
\label{3.2}
The smooth solutions of system \ref{3systeme2} satisfy the following inequality when $P(\rho)=a\rho^{\gamma}$ with $\gamma\geq 1$:
\begin{equation}
\begin{aligned}
&\frac{d}{dt}\int_{\Om}\rho\frac{|v|^{2+\de}}{2+\de}+\frac{\nu}{4}\int_{\Om} \rho|v|^{\de}|\n v|^{2}dx\leq
\biggl(\int_{\Om}\big(\rho^{2\gamma-1-\frac{\delta}{2}}\big)^{\frac{2}{2-\de}}dx\biggl)^{\frac{2}{2-\de}}
\big(\int_{\Om}\rho|v|^{2}dx\big)^{\frac{\de}{2}},
\end{aligned}
\label{22a}
\end{equation}
for $\delta\in(0,2)$.
\end{proposition}
\section{Littlewood-Paley theory and Besov spaces}
Throughout the paper, $C$ stands for a constant whose exact meaning depends on the context. The notation $A\lesssim B$  means
that $A\leq CB$.
For all Banach space $X$, we denote by $C([0,T],X)$ the set of continuous functions on $[0,T]$ with values in $X$.
For $p\in[1,+\infty]$, the notation $L^{p}(0,T,X)$ or $L^{p}_{T}(X)$ stands for the set of measurable functions on $(0,T)$
with values in $X$ such that $t\rightarrow\|f(t)\|_{X}$ belongs to $L^{p}(0,T)$.
Littlewood-Paley decomposition  corresponds to a dyadic
decomposition  of the space in Fourier variables.
Let $\varphi\in C^{\infty}(\mathbb{R}^{N})$,
supported in the shell
${\cal{C}}=\{\xi\in\R^{N}/\frac{3}{4}\leq|\xi|\leq\frac{8}{3}\}$ and  $\chi\in C^{\infty}(\mathbb{R}^{N})$
supported in the ball $B(0,\frac{4}{3})$. $\varphi$ and  $\chi$ are valued in $[0,1]$.
We set ${\cal Q}^{N}=(0,2\pi)^{N}$ and $\widetilde{\mathbb{Z}}^{N}=(\mathbb{Z}/1)^{N}$ the dual lattice associated to $\mathbb{T}^{N}$. We decompose now $u\in{\cal S}^{'}(\T^{N})$ into Fourier series:
$$u(x)=\sum_{\beta\in \widetilde{\mathbb{Z}}^{N}}\hat{u}_{\beta}e^{i\beta\cdot x}\;\;\;\mbox{with}\;\;
\hat{u}_{\beta}=\frac{1}{|\T^{N}|}\int_{\T^{N}}e^{-i\beta\cdot y}u(y)dy.$$
Denoting;
$$h_{q}(x)=\sum_{\beta\in \widetilde{\mathbb{Z}}^{N}}\va(2^{-q}\beta)e^{i\beta\cdot x},$$
one can now define the periodic dyadic blocks as:
$$\D_{q}u(x)=\sum_{\beta\in \widetilde{\mathbb{Z}}^{N}}\va(2^{-q}\beta)\hat{u}_{\beta}e^{i\beta\cdot x}=\frac{1}{|\T^{N}|}\int_{\T^{N}}h_{q}(y)u(x-y)dy,\;\;\mbox{for all}\;\;q\in\mathbb{Z}$$
and the low frequency cutt-off:
$$S_{q}u(x)=\hat{u}_{0}+\sum_{p\leq q-1}\D_{p}u(x)=\sum_{\beta\in \widetilde{\mathbb{Z}}^{N}}\chi(2^{-q}\beta)\hat{u}_{\beta}e^{i\beta\cdot x}.$$
It is obvious that:
$$u=\hat{u}_{0}+\sum_{k\in\mathbb{Z}}\D_{k}u.$$
This decomposition is called non-homogeneous Littlewood-Paley
decomposition.
\\
Furthermore we have the following proposition where $\widetilde{{\cal C}}=B(0,\frac{2}{3})+{\cal{C}}$
\begin{proposition}
\label{d210}
\begin{equation}
|k-k^{'}|\geq 2\implies {\rm supp}\va(2^{-k}\cdot)\cap {\rm supp}\va(2^{-k^{'}}\cdot)=\emptyset,
\label{d28a}
\end{equation}
\begin{equation}
k\geq 1\implies {\rm supp}\chi\cap {\rm supp}\va(2^{-k}\cdot)=\emptyset,
\label{d29a}
\end{equation}
\begin{equation}
|k-k^{'}|\geq 5\implies 2^{k^{'}}\widetilde{{\cal C}}\cap 2^{k}{\cal{C}}=\emptyset.
\label{d210a}
\end{equation}
\end{proposition}
\subsection{Non homogeneous Besov spaces and first properties}
\begin{definition}
For
$s\in\R,\,\,p\in[1,+\infty],\,\,q\in[1,+\infty],\,\,\mbox{and}\,\,u\in{\cal{S}}^{'}(\R^{N})$
we set:
$$\|u\|_{B^{s}_{p,q}}=(\sum_{l\in\mathbb{Z}}(2^{ls}\|D_{l}u\|_{L^{p}})^{q})^{\frac{1}{q}}.$$
The non homogeneous Besov space $B^{s}_{p,q}$ is the set of temperate  distribution $u$  such that $\|u\|_{\dot{B}^{s}_{p,q}}<+\infty$.
\end{definition}
\begin{remarka}The above definition is a natural generalization of the
nonhomogeneous Sobolev and H$\ddot{\mbox{o}}$lder spaces: one can show
that $B^{s}_{\infty,\infty}$ is the nonhomogeneous
H$\ddot{\mbox{o}}$lder space $C^{s}$ and that $B^{s}_{2,2}$ is
the nonhomogeneous space $H^{s}$.
\end{remarka}
\begin{proposition}
\label{derivation,interpolation}
The following properties holds:
\begin{enumerate}
\item there exists a constant universal $C$
such that:\\
$C^{-1}\|u\|_{B^{s}_{p,r}}\leq\|\n u\|_{B^{s-1}_{p,r}}\leq
C\|u\|_{B^{s}_{p,r}}.$
\item If
$p_{1}<p_{2}$ and $r_{1}\leq r_{2}$ then $B^{s}_{p_{1},r_{1}}\hookrightarrow
B^{s-N(1/p_{1}-1/p_{2})}_{p_{2},r_{2}}$.
\item $B^{s^{'}}_{p,r_{1}}\hookrightarrow B^{s}_{p,r}$ if $s^{'}> s$ or if $s=s^{'}$ and $r_{1}\leq r$.
\end{enumerate}
\label{interpolation}
\end{proposition}
Let now recall a few product laws in Besov spaces coming directly from the paradifferential calculus of J-M. Bony
(see \cite{fBJM}) and rewrite on a generalized form in \cite{AP} by H. Abidi and M. Paicu (in this article the results are written
in the case of homogeneous sapces but it can easily generalize for the nonhomogeneous Besov spaces).
\begin{proposition}
\label{produit1}
We have the following laws of product:
\begin{itemize}
\item For all $s\in\R$, $(p,r)\in[1,+\infty]^{2}$ we have:
\begin{equation}
\|uv\|_{B^{s}_{p,r}}\leq
C(\|u\|_{L^{\infty}}\|v\|_{B^{s}_{p,r}}+\|v\|_{L^{\infty}}\|u\|_{B^{s}_{p,r}})\,.
\label{a2.2}
\end{equation}
\item Let $(p,p_{1},p_{2},r,\lambda_{1},\lambda_{2})\in[1,+\infty]^{2}$ such that:$\frac{1}{p}\leq\frac{1}{p_{1}}+\frac{1}{p_{2}}$,
$p_{1}\leq\lambda_{2}$, $p_{2}\leq\lambda_{1}$, $\frac{1}{p}\leq\frac{1}{p_{1}}+\frac{1}{\lambda_{1}}$ and
$\frac{1}{p}\leq\frac{1}{p_{2}}+\frac{1}{\lambda_{2}}$. We have then the following inequalities:\\
if $s_{1}+s_{2}+N\inf(0,1-\frac{1}{p_{1}}-\frac{1}{p_{2}})>0$, $s_{1}+\frac{N}{\lambda_{2}}<\frac{N}{p_{1}}$ and
$s_{2}+\frac{N}{\lambda_{1}}<\frac{N}{p_{2}}$ then:
\begin{equation}
\|uv\|_{B^{s_{1}+s_{2}-N(\frac{1}{p_{1}}+\frac{1}{p_{2}}-\frac{1}{p})}_{p,r}}\lesssim\|u\|_{B^{s_{1}}_{p_{1},r}}
\|v\|_{B^{s_{2}}_{p_{2},\infty}},
\label{a2.3}
\end{equation}
when $s_{1}+\frac{N}{\lambda_{2}}=\frac{N}{p_{1}}$ (resp $s_{2}+\frac{N}{\lambda_{1}}=\frac{N}{p_{2}}$) we replace
$\|u\|_{B^{s_{1}}_{p_{1},r}}\|v\|_{B^{s_{2}}_{p_{2},\infty}}$ (resp $\|v\|_{B^{s_{2}}_{p_{2},\infty}}$) by
$\|u\|_{B^{s_{1}}_{p_{1},1}}\|v\|_{B^{s_{2}}_{p_{2},r}}$ (resp $\|v\|_{B^{s_{2}}_{p_{2},\infty}\cap L^{\infty}}$),
if $s_{1}+\frac{N}{\lambda_{2}}=\frac{N}{p_{1}}$ and $s_{2}+\frac{N}{\lambda_{1}}=\frac{N}{p_{2}}$ we take $r=1$.
\\
If $s_{1}+s_{2}=0$, $s_{1}\in(\frac{N}{\lambda_{1}}-\frac{N}{p_{2}},\frac{N}{p_{1}}-\frac{N}{\lambda_{2}}]$ and
$\frac{1}{p_{1}}+\frac{1}{p_{2}}\leq 1$ then:
\begin{equation}
\|uv\|_{B^{-N(\frac{1}{p_{1}}+\frac{1}{p_{2}}-\frac{1}{p})}_{p,\infty}}\lesssim\|u\|_{B^{s_{1}}_{p_{1},1}}
\|v\|_{B^{s_{2}}_{p_{2},\infty}}.
\label{a2.4}
\end{equation}
If $|s|<\NN$ for $p\geq2$ and $-\frac{N}{p^{'}}<s<\NN$ else, we have:
\begin{equation}
\|uv\|_{B^{s}_{p,r}}\leq C\|u\|_{B^{s}_{p,r}}\|v\|_{B^{\NN}_{p,\infty}\cap L^{\infty}}.
\label{a2.5}
\end{equation}
\end{itemize}
\end{proposition}
\begin{remarka}
In the sequel $p$ will be either $p_{1}$ or $p_{2}$ and in this case $\frac{1}{\lambda}=\frac{1}{p_{1}}-\frac{1}{p_{2}}$
if $p_{1}\leq p_{2}$, resp $\frac{1}{\lambda}=\frac{1}{p_{2}}-\frac{1}{p_{1}}$
if $p_{2}\leq p_{1}$.
\end{remarka}
\begin{corollaire}
\label{produit2}
Let $r\in [1,+\infty]$, $1\leq p\leq p_{1}\leq +\infty$ and $s$ such that:
\begin{itemize}
\item $s\in(-\frac{N}{p_{1}},\frac{N}{p_{1}})$ if $\frac{1}{p}+\frac{1}{p_{1}}\leq 1$,
\item $s\in(-\frac{N}{p_{1}}+N(\frac{1}{p}+\frac{1}{p_{1}}-1),\frac{N}{p_{1}})$ if $\frac{1}{p}+\frac{1}{p_{1}}> 1$,
\end{itemize}
then we have if $u\in B^{s}_{p,r}$ and $v\in B^{\frac{N}{p_{1}}}_{p_{1},\infty}\cap L^{\infty}$:
$$\|uv\|_{B^{s}_{p,r}}\leq C\|u\|_{B^{s}_{p,r}}\|v\|_{B^{\frac{N}{p_{1}}}_{p_{1},\infty}\cap L^{\infty}}.$$
\end{corollaire}
The study of non stationary PDE's requires space of type $L^{\rho}(0,T,X)$ for appropriate Banach spaces $X$. In our case, we
expect $X$ to be a Besov space, so that it is natural to localize the equation through Littlewood-Payley decomposition. But, in doing so, we obtain
bounds in spaces which are not type $L^{\rho}(0,T,X)$ (except if $r=p$).
We are now going to
define the spaces of Chemin-Lerner in which we will work, which are
a refinement of the spaces
$L_{T}^{\rho}(B^{s}_{p,r})$.
$\hspace{15cm}$
\begin{definition}
Let $\rho\in[1,+\infty]$, $T\in[1,+\infty]$ and $s_{1}\in\R$. We set:
$$\|u\|_{\widetilde{L}^{\rho}_{T}(B^{s_{1}}_{p,r})}=
\big(\sum_{l\in\mathbb{Z}}2^{lrs_{1}}\|\D_{l}u(t)\|_{L^{\rho}(L^{p})}^{r}\big)^{\frac{1}{r}}\,.$$
We then define the space $\widetilde{L}^{\rho}_{T}(B^{s_{1}}_{p,r})$ as the set of temperate distribution $u$ over
$(0,T)\times\R^{N}$ such that 
$\|u\|_{\widetilde{L}^{\rho}_{T}(B^{s_{1}}_{p,r})}<+\infty$.
\end{definition}
We set $\widetilde{C}_{T}(\widetilde{B}^{s_{1}}_{p,r})=\widetilde{L}^{\infty}_{T}(\widetilde{B}^{s_{1}}_{p,r})\cap
{\cal C}([0,T],B^{s_{1}}_{p,r})$.
Let us emphasize that, according to Minkowski inequality, we have:
$$\|u\|_{\widetilde{L}^{\rho}_{T}(B^{s_{1}}_{p,r})}\leq\|u\|_{L^{\rho}_{T}(B^{s_{1}}_{p,r})}\;\;\mbox{if}\;\;r\geq\rho
,\;\;\;\|u\|_{\widetilde{L}^{\rho}_{T}(B^{s_{1}}_{p,r})}\geq\|u\|_{L^{\rho}_{T}(B^{s_{1}}_{p,r})}\;\;\mbox{if}\;\;r\leq\rho
.$$
\begin{remarka}
It is easy to generalize proposition \ref{produit1},
to $\widetilde{L}^{\rho}_{T}(B^{s_{1}}_{p,r})$ spaces. The indices $s_{1}$, $p$, $r$
behave just as in the stationary case whereas the time exponent $\rho$ behaves according to H\"older inequality.
\end{remarka}
In the sequel we will need of composition lemma in $\widetilde{L}^{\rho}_{T}(B^{s}_{p,r})$ spaces.
\begin{lemme}
\label{composition}
Let $s>0$, $(p,r)\in[1,+\infty]$ and $u\in \widetilde{L}^{\rho}_{T}(B^{s}_{p,r})\cap L^{\infty}_{T}(L^{\infty})$.
\begin{enumerate}
 \item Let $F\in W_{loc}^{[s]+2,\infty}(\R^{N})$ such that $F(0)=0$. Then $F(u)\in \widetilde{L}^{\rho}_{T}(B^{s}_{p,r})$. More precisely there exists a function $C$ depending only on $s$, $p$, $r$, $N$ and $F$ such that:
$$\|F(u)\|_{\widetilde{L}^{\rho}_{T}(B^{s}_{p,r})}\leq C(\|u\|_{L^{\infty}_{T}(L^{\infty})})\|u\|_{\widetilde{L}^{\rho}_{T}(B^{s}_{p,r})}.$$
\item Let $F\in W_{loc}^{[s]+3,\infty}(\R^{N})$ such that $F(0)=0$. Then $F(u)-F^{'}(0)u\in \widetilde{L}^{\rho}_{T}(B^{s}_{p,r})$. More precisely there exists a function $C$ depending only on $s$, $p$, $r$, $N$ and $F$ such that:
$$\|F(u)-F^{'}(0)u\|_{\widetilde{L}^{\rho}_{T}(B^{s}_{p,r})}\leq C(\|u\|_{L^{\infty}_{T}(L^{\infty})})\|u\|^{2}_{\widetilde{L}^{\rho}_{T}(B^{s}_{p,r})}.$$
\end{enumerate}
\end{lemme}
Let us now give some estimates for the heat equation:
\begin{proposition}
\label{5chaleur} Let $s\in\R$, $(p,r)\in[1,+\infty]^{2}$ and
$1\leq\rho_{2}\leq\rho_{1}\leq+\infty$. Assume that $u_{0}\in B^{s}_{p,r}$ and $f\in\widetilde{L}^{\rho_{2}}_{T}
(B^{s-2+2/\rho_{2}}_{p,r})$.
Let u be a solution of:
$$
\begin{cases}
\begin{aligned}
&\p_{t}u-\mu\D u=f\\
&u_{t=0}=u_{0}\,.
\end{aligned}
\end{cases}
$$
Then there exists $C>0$ depending only on $N,\mu,\rho_{1}$ and
$\rho_{2}$ such that:
$$\|u\|_{\widetilde{L}^{\rho_{1}}_{T}(\widetilde{B}^{s+2/\rho_{1}}_{p,r})}\leq C\big(
 \|u_{0}\|_{B^{s}_{p,r}}+\mu^{\frac{1}{\rho_{2}}-1}\|f\|_{\widetilde{L}^{\rho_{2}}_{T}
 (B^{s-2+2/\rho_{2}}_{p,r})}\big)\,.$$
 If in addition $r$ is finite then $u$ belongs to $C([0,T],B^{s}_{p,r})$.
\end{proposition}
\section{Proof of  the theorems \ref{2.1} and \ref{2.2}}
We now present the proof of theorem \ref{2.1}. To begin with, we need to make precise the assumptions on the initial data.
\subsubsection*{Initial data:}
We recall that the initial data must satisfy (\ref{16}), and (\ref{17}) to make use of all the inequalities presented in the previous section:
\begin{itemize}
\item $\rho_{0}^{n}$ is bounded in $L^{1}(\Om)\cap L^{\gamma}(\Om)$, $\rho_{0}^{n}\geq0$ a.e in $\Om$,
\item $\rho_{0}^{n}|u_{0}^{n}|^{2}$ is bounded in $L^{1}(\Om)$,
\item $\sqrt{\rho_{0}^{n}}\n\ln\rho_{0}^{n}$ is bounded in $L^{2}(\Om)$,
\item $\rho_{0}^{n}|v_{0}^{n}|^{2+\delta}$ is bounded in $L^{1}(\Om)$.
\end{itemize}
With those assumptions, and using the entropies (\ref{21}), (\ref{22}) and the mass equation, we deduce the following energy estimates, which we shall give uniform bounds:
\begin{equation}
\begin{aligned}
&\|\sqrt{\rho_{n}}\|_{L^{\infty}((0,T),L^{2}(\Om))}\leq C,\\
&\|\rho_{n}\|_{L^{\infty}((0,T),L^{\gamma}(\Om))}\leq C,\\
&\|\sqrt{\rho_{n}}\p_{ij}\ln\rho_{n}\|_{L^{2}((0,T)\times \Om)}\leq C,\\
&\|\sqrt{\rho_{n}}\n u_{n}\|_{L^{2}((0,T)\times\Om}\leq C,
\end{aligned}
\label{27}
\end{equation}
and
\begin{equation}
\begin{aligned}
&\|\rho_{n}^{\frac{\gamma}{2}-1}\n\rho_{n}\|_{L^{2}((0,T)\times\Om)}\leq C,\\
&\|\rho_{n}|v_{n}|^{2+\delta}\|_{L^{\infty}((0,T),L^{1}(\Om))} \leq C.
\end{aligned}
\label{28}
\end{equation}
In view of previous inequalities, the bounds (\ref{27}) and (\ref{28}) yields:
\begin{equation}
\begin{aligned}
&\|\sqrt{\rho_{n}}\n v_{n}\|_{L^{2}((0,T)\times\Om)}\leq C,\\
&\|\n\sqrt{\rho_{n}}\|_{L^{\infty}(0,T;L^{2}(\Om)}\leq C,\\
&\|\n\rho_{n}^{\frac{\gamma}{2}}\|_{L^{2}((0,T)\times\Om)}\leq C,\\
&\|\rho_{n}|v_{n}|^{2+\delta}\|_{L^{\infty}((0,T),L^{1}(\Om))} \leq C.
\end{aligned}
\label{29}
\end{equation}
The proof of theorem \ref{2.1} will be derived in $5$ steps and follows the proof of \cite{fMV}. In the first two steps, we show the convergence of the density and the pressure (note that the convergence of the density is straightforward here). The key argument of the proof is presented in the third step: as $\sqrt{\rho_{n}}v_{n}$ is bounded in a space better than $L^{\infty}(0,T;L^{2}(\R))$, it will allow to give the convergence of the momentum (step 4) and finally the strong convergence of $\sqrt{\rho_{n}}v_{n}\otimes \sqrt{\rho_{n}}u_{n}$ in $L_{w}^{2}((0,T)\times\R)$ (step 5). In fact by considering the effective velocity $v_{n}$ we do not need extra information on the density $\rho_{n}$ to deal with the capillarity terms (as it is the case in \cite{5BDL,fH2}). Furthermore at the difference with \cite{5BDL} we are able to treat the momentum term without using test functions depending of the density (or in an other way to add a cold pressure).
The last step address the convergence of the diffusion terms and of $\sqrt{\rho_{n}}v_{n}$ to
$\sqrt{\rho}u+\frac{2\kappa}{\mu}\n\rho$.
\subsubsection*{Step 1: Convergence of $\sqrt{\rho_{n}}$}
\begin{lemme}
We have the following properties:
\begin{enumerate}
\item $\sqrt{\rho_{n}}$ is bounded in $L^{\infty}(0,T;H^{1}(\Om))$
\item $\p_{t}\sqrt{\rho_{n}}$ is bounded in $L^{2}(0,T;H^{-1}(\Om))$.
\end{enumerate}
As a consequence up to a subsequence, $\sqrt{\rho_{n}}$ converge almost everywhere and strongly in $L^{2}(0,T;L^{2}_{loc}(\Om))$. We write:
$$\sqrt{\rho_{n}}\rightarrow\sqrt{\rho}\;\;\;\mbox{a.e and}\;\;L^{2}_{loc}((0,T)\times \Om)\;\;\mbox{strong}.$$
Moreover, $\rho_{n}$ converges to $\rho$ in $C^{0}(0,T;L^{\frac{3}{2}}_{loc}(\Om))$.
\label{lemme1}
\end{lemme}
{\bf Proof:} The second estimate in (\ref{29}), together with the conservation of mass $\|\rho_{n}\|_{L^{1}}=\|\rho_{n,0}\|_{L^{1}}$ gives the $L^{\infty}(0,T;H^{1}(\Om))$
bound. 
Next as we have:
$$
\begin{aligned}
\p_{t}\sqrt{\rho_{n}}&=-\frac{1}{2}\sqrt{\rho_{n}}{\rm div}u_{n}-u_{n}\cdot\n\sqrt{\rho_{n}}\\
&=\frac{1}{2}\sqrt{\rho_{n}}{\rm div}u_{n}-{\rm div}(u_{n}\sqrt{\rho_{n}}),
\end{aligned}
$$
 thanks to Aubin's Lions Lemma,  we conclude to the strong convergence in $L^{2}_{loc}((0,T)\times\Om)$ of $\sqrt{\rho_{n}}$. Up to a subsequence we obtain that $\rho_{n}$ converges to $\rho$ a.e.\\
Sobolev embedding implies that $\sqrt{\rho_{n}}$ is bounded in $L^{\infty}(0,T;L^{q}(\Om))$ for $q\in[2,+\infty[$ if $N=2$ and $q\in[2,6]$ if $N=3$.
In either cases we deduce that $\rho_{n}$ is bounded in $L^{\infty}(0,T;L^{3}(\Om))$, and therefore:
$$\rho_{n}u_{n}=\sqrt{\rho_{n}}\sqrt{\rho_{n}}u_{n}\;\;\mbox{is bounded in}\;\;L^{\infty}(0,T;L^{\frac{3}{2}}(\Om)).$$
The continuity equation thus yields $\p_{t}\rho_{n}$ bounded in $L^{\infty}(0,T;W^{-1,\frac{3}{2}}(\Om))$. Moreover since
$\n\rho_{n}=2\sqrt{\rho_{n}}\n \sqrt{\rho_{n}}$, we also have that $\n\rho_{n}$, is bounded in $L^{\infty}(0,T;L^{\frac{3}{2}}(\Om))$, hence the compactness of $\rho_{n}$ in $C([0,T];L^{\frac{3}{2}}_{loc}(\Om))$.
\subsubsection*{Step 2: Convergence of the pressure}
Our goal now is to prove that the pressure term is enough integrable in $L^{1+\e}((0,T)\times\Omega)$. Indeed it will be sufficient to prove that $P(\rho_{n})$ converges to $P(\rho)$ in $L^{1}_{loc}(0,T)\times\Omega)$ since we know from the previous lemma that $\rho_{n}$ converges a.e to $\rho$.
\begin{lemme}
The pressure $\rho_{n}^{\gamma}$ is bounded in $L^{\frac{5}{3}}((0,T)\times\Omega)$ when $N=3$ and $L^{r}((0,T)\times\Omega)$
for all $r\in[1,2[$ when $N=2$. In particular, $\rho_{n}^{\gamma}$ converges to $\rho^{\gamma}$ strongly in $L^{1}_{loc}((0,T)\times\Om)$.
\label{pression}
\end{lemme}
{\bf Proof:} Inequalities (\ref{29}) and (\ref{28}) yield $\rho_{n}^{\frac{\gamma}{2}}\in L^{2}(0,T;H^{1})$.\\
When $N=2$, we deduce $\rho_{n}^{\frac{\gamma}{2}}\in L^{2}(0,T;L^{q}(\R^{N}))$ for all $q\in[2,\infty[$. So $\rho_{n}^{\gamma}$
is bounded in $L^{1}(0,T;L^{p})\cap L^{\infty}(L^{1})$ for all $p\in[1,+\infty[$, hence by interpolation $\rho_{n}^{\gamma}$
is bounded in $L^{r}((0,T)\times\R)$ for all $r\in[1,2[$.\\
When $N=3$, we  get by Sobolev embedding that $\rho_{n}^{\gamma}$ is bounded in $L^{1}(0,T;L^{3}(\R))$. As $\rho_{n}^{\gamma}$ is also in $L^{\infty}((0,T),L^{1}(\Om))$, by interpolation we have:
$$\|\rho_{n}^{\gamma}\|_{L^{\frac{5}{3}}((0,T)\times\Om)}\leq\|\rho_{n}^{\gamma}\|_{L^{\infty}(0,T;L^{1}(\Om))}^{\frac{2}{5}}
\|\rho_{n}^{\gamma}\|_{L^{1}(0,T;L^{3}(\Om))}^{\frac{3}{5}}$$
Hence $\rho_{n}^{\gamma}$ is bounded in $L^{\frac{5}{3}}((0,T)\times\Om)$.\\
\\
Since we already now that $\rho_{n}^{\gamma}$ converges almost everywhere to $\rho^{\gamma}$, we can conclude to the strong convergence of $\rho_{n}^{\gamma}$ in $L^{1}_{loc}((0,T)\times\R)$.
\subsubsection*{Step 3: Convergence of the momentum}
\begin{lemme}
Up to a subsequence, the momentum $m_{n}=\rho_{n}v_{n}$ converges strongly in $L^{2}(0,T;L^{p}_{loc}(\Om))$ to some $m(x,t)$ for
all $p\in[1,\frac{3}{2})$. In particular:
$$\rho_{n}v_{n}\rightarrow m\;\;\;\mbox{almost everywhere}\;\;(x,t)\in(0,T)\times\Om.$$
Note that we can already define $v(t,x)=\frac{m(t,x)}{\rho(t,x)}$ outside the vacuum set $\{\rho(t,x)=0\}$,
but we do not know yet whether $m(t,x)$ is zero on the vacuum set.
\label{moment}
\end{lemme}
{\bf Proof:} We have:
$$\rho_{n}v_{n}=\sqrt{\rho_{n}}\sqrt{\rho_{n}}v_{n},$$
where $\sqrt{\rho_{n}}$ is bounded in $L^{\infty}(0,T;L^{q}(\Om))$ for $q\in[2,+\infty[$ if $N=2$ and $q\in[2,6]$ if $N=3$.
Since $\sqrt{\rho_{n}}v_{n}$ is bounded in $L^{\infty}(0,T;L^{2}(\Om))$, we deduce that:
$$\rho_{n}v_{n}\;\;\mbox{is bounded in}\;\;L^{\infty}(0,T;L^{q}(\Om))\;\;\mbox{for all}\;\;q\in[1,\frac{3}{2}].$$
Next we have:
$$
\begin{aligned}
\p_{i}(\rho_{n}v^{j}_{n})&=\rho_{n}\p_{i}v_{nj}+u^{j}_{n}\p_{i}\rho_{n}\\
&=\sqrt{\rho_{n}}\sqrt{\rho_{n}}\p_{i}v^{j}_{n}+2\sqrt{\rho_{n}}v_{n}\p_{i}\sqrt{\rho_{n}}.\\
\end{aligned}
$$
Using (\ref{29}) the second term is bounded in $L^{\infty}(0,T;L^{1}(\Om))$, while the first
term is bounded in $L^{2}(0,T;L^{q}(\Om))$ for all $q\in[1,\frac{3}{2}]$. Then:
$$\n(\rho_{n}v_{n})\;\;\mbox{is bounded in}\;\;L^{2}(0,T;L^{1}(\Om)).$$
In particular we have:
$$\rho_{n}v_{n}\;\;\mbox{is bounded in}\;\;L^{2}(0,T;W^{1,1}(\Om)).$$
It remains to show that for every compact set $K\in\Om$, we have:
\begin{equation}
\p_{t}(\rho_{n}v_{n})\;\;\mbox{is bounded in}\;\;L^{2}(0,T;W^{-2,\frac{4}{3}}(K)).
\label{32}
\end{equation}
Indeed (\ref{32}) and Aubin's Lemma (see Lions \cite{Li}), yields the compactness of $\rho_{n}v_{n}$ in $L^{2}(0,T;L^{p}(K))$ for all $p\in[1,\frac{3}{2})$.\\
To prove (\ref{32}), we use the momentum equation of system (\ref{3systeme1}), first we observe that using  lemma \ref{pression} and (\ref{29}):
$$
\begin{aligned}
{\rm div}(\sqrt{\rho_{n}}v_{n}\otimes\sqrt{\rho_{n}}u_{n})&\in L^{\infty}(0,T;W^{-1,1}(K)),\\
\n\rho_{n}^{\gamma}&\in L^{\infty}(0,T;W^{-1,1}(K)).
\end{aligned}
$$
So we only have to check that the terms $\n(\rho_{n}\n v_{n})$, $\n(\rho_{n}^{t}\n v_{n})$ are bounded in
$L^{\infty}(0,T;W^{-2,\frac{4}{3}}(K))$. To that purpose, we write:
\begin{equation}
\rho_{n}\n v_{n}=\n(\rho_{n}v_{n})-v_{n}\n\rho_{n},
\label{33}
\end{equation}
(and we proceed similarly with the other two terms). The second term in (\ref{33}) reads:
$$v_{n}\n\rho_{n}=2\sqrt{\rho_{n}}v_{n}\n\sqrt{\rho_{n}},$$
which is bounded in $L^{\infty}(0,T;L^{1}(\Omega))$ thanks to (\ref{28}) and (\ref{29}). The first term in (\ref{33}) can be rewritten:
$$\n[\rho_{n}v_{n}]=\n[\sqrt{\rho_{n}}(\sqrt{\rho_{n}}v_{n})],$$
which is bounded in $L^{\infty}(0,T;W^{-1,\frac{3}{2}}(\Omega))$ thanks to lemma \ref{lemme1}.\\
We deduce that $\rho_{n}D(v_{n})$ is bounded in:
$$L^{\infty}(0,T;W^{-1,\frac{3}{2}}(K)+L^{1}(K)),$$
and since $L^{1}(K)\subset W^{-1,\frac{4}{3}}(K)$ and $W^{-1,\frac{3}{2}}(K)\subset W^{-1,\frac{4}{3}}(K)$ we conclude that $\rho_{n}\,D(v_{n})$ and $\rho_{n}\,\n v_{n}$ are bounded in $L^{\infty}(0,T;W^{-1,\frac{4}{3}}(K))$, which conclude the proof of lemma \ref{moment}.\\
\subsubsection*{Step 4: Convergence of $\sqrt{\rho_{n}}u_{n}\otimes\sqrt{\rho_{n}}v_{n}$}
\begin{lemme}
The quantity $\sqrt{\rho_{n}}v_{n}$converges strongly in $L^{2}_{loc}((0,T)\times\Omega)$ to $\frac{m}{\sqrt{\rho}}$
(defined to be zero when $m=0$).\\
In particular, we have $m(t,x)=0$ a.e on $\{\rho(t,x)=0\}$ and there exists a function $v(t,x)$ such that $m(t,x)=\rho(t,x)v(t,x)$ and:
$$
\begin{aligned}
&\rho_{n}v_{n}\rightarrow\rho v\;\;\;\mbox{strongly in}\;\;L^{2}(0,T;L^{p}_{loc}(\Omega)),\;p\in[1,\frac{3}{2}),\\
&\sqrt{\rho_{n}}v_{n}\rightarrow\sqrt{\rho}v\;\;\;\mbox{strongly in}\;\;L^{2}_{loc}((0,T)\times\Omega),\\
\end{aligned}
$$
(note that $v$ is not uniquely defined on the vacuum set $\{\rho(t,x)=0\}$).
\label{imp1}
\end{lemme}
{\bf Proof:} First of all, since $\frac{m_{n}}{\sqrt{\rho_{n}}}$ is bounded in $L^{\infty}(0,T;L^{2}(\Omega))$, Fatou's lemma yields:
$$\int\lim\inf\frac{m_{n}^{2}}{\rho_{n}}dx<+\infty.$$
In particular, we have $m(t,x)=0$ a.e. in $\{\rho(t,x)=0\}$. So if we define the limit velocity by $v(t,x)$ by setting $v(t,x)=\frac{m(t,x)}{\rho(t,x)}$ when $\rho(t,x)\ne0$ and $v(t,x)=0$ when $\rho(t,x)=0$, we have:
$$m(t,x)=\rho(t,x)v(t,x)$$
and
$$\int\frac{m^{2}}{\rho}dx=\int\rho|v|^{2}dx<+\infty.$$
Moreover, Fatou's lemma yields:
$$
\begin{aligned}
\int\rho|v|^{2+\delta}dx&\leq\int\lim\inf\rho_{n}|v_{n}|^{2+\delta}dx\\
&\leq\lim\inf\int\rho_{n}|v_{n}|^{2+\delta}dx,
\end{aligned}
$$
and so $\rho|v|^{2+\delta}$ is in $L^{\infty}(0,T;L^{1}(\Omega))$.\\
Next, since $m_{n}$ and $\rho_{n}$ converge almost everywhere, it is readily seen that in $\{\rho(t,x)\ne0\}$, $\sqrt{\rho_{n}}v_{n}=\frac{m_{n}}
{\sqrt{\rho_{n}}}$ converges almost everywhere to $\sqrt{\rho}u=\frac{m}{\sqrt{\rho}}$. Moreover, we have:
\begin{equation}
\sqrt{\rho_{n}}v_{n}1_{\{|v_{n}|\leq M\}}\rightarrow\sqrt{\rho}v1_{\{|v|\leq M\}}\;\;\;\mbox{almost everywhere}.
\label{34}
\end{equation}
As a matter of fact, the convergence holds almost everywhere in $\{\rho(t,x)\ne0\}$, and in $\{\rho(t,x)=0\}$, we have
$\sqrt{\rho_{n}}v_{n}1_{\{|v_{n}|\leq M\}}\leq M\sqrt{\rho_{n}}\rightarrow0$.\\
We are now in position to complete the proof of lemma \ref{imp1}. For $M>0$ we cut the $L^{2}$ norm as follows:
$$
\begin{aligned}
&\int|\sqrt{\rho_{n}}v_{n}-\sqrt{\rho}v|^{2}dxdt\leq\int|\sqrt{\rho_{n}}v_{n}1_{\{|v_{n}|\leq M\}}-\sqrt{\rho}v1_{\{|v|\leq M\}}|^{2}dxdt\\[2mm]
&\hspace{3cm}+2\int|\sqrt{\rho_{n}}v_{n}1_{\{|v_{n}|\geq M\}}|^{2}dxdt+2\int|\sqrt{\rho}v1_{\{|v|\geq M\}}|^{2}dxdt,
\end{aligned}
$$
It is obvious that $\sqrt{\rho_{n}}v_{n}1_{\{|v_{n}|\leq M\}}$ is bounded uniformly in $L^{\infty}(0,T;L^{3}(\Omega))$ and by (\ref{34}), we obtain that $\sqrt{\rho_{n}}v_{n}1_{\{|v_{n}|\leq M\}}$ converges strongly to $\sqrt{\rho}v1_{\{|v||\leq M\}}$ in $L^{2}_{loc}((0,T)\times\Om)$:
\begin{equation}
\int|\sqrt{\rho_{n}}v_{n}1_{\{|v_{n}|\leq M\}}-\sqrt{\rho}v1_{\{|v|\leq M\}}|^{2}dxdt\rightarrow0.
\label{35}
\end{equation}
Finally, we write:
\begin{equation}
\int|\sqrt{\rho_{n}}v_{n}1_{\{|v_{n}|\geq M\}}|^{2}dxdt\leq\frac{1}{M^{\delta}}\int\rho|v_{n}|^{2+\delta}dxdt,
\label{36}
\end{equation}
and
\begin{equation}
\int|\sqrt{\rho}v1_{\{|v|\geq M\}}|^{2}dxdt\leq\frac{1}{M^{\delta}}\int\rho|v|^{2+\delta}dxdt,
\label{37}
\end{equation}
Putting together (\ref{35}), (\ref{36}) and (\ref{37}), we deduce:
$$\lim\sup_{n\rightarrow+\infty}\int|\sqrt{\rho_{n}}u_{n}-\sqrt{\rho}u|^{2}dxdt\leq\frac{C}{M^{\delta}},$$
for all $M>0$, and the lemma follows by taking $M\rightarrow+\infty$.
\begin{lemme}
$\sqrt{\rho_{n}}u_{n}\otimes\sqrt{\rho_{n}}v_{n}$ converges to the distribution sense to
$\sqrt{\rho_{n}}u_{n}\otimes\sqrt{\rho_{n}}v_{n}$.
\end{lemme}
{\bf Proof:} It suffices to observe that $\sqrt{\rho_{n}}u_{n}$ converges weakly in $L^{2}(0,T,L^{2}(\Om))$ to $\sqrt{\rho}u$ and $\sqrt{\rho_{n}}v_{n}$ converges strongly in $L^{2}_{loc}((0,T)\times \Om)$ to $\sqrt{\rho}v$. We have then for all $\va \in C^{\infty}_{0}((0,T)\times\Om)$:
$$\int^{T}_{0}\int_{\Om}(\sqrt{\rho_{n}}u_{n}\otimes \sqrt{\rho_{n}}v_{n})\va dxdt\rightarrow
\int^{T}_{0}\int_{\Om}(\sqrt{\rho}u\otimes \sqrt{\rho}v)\va dxdt.$$
\subsubsection*{ Step 5: Convergence of the diffusion terms}
\begin{lemme}
We have:
$$
\begin{aligned}
&\rho_{n}\n v_{n}\rightarrow \rho\n v\;\;\;\mbox{in}\;\;{\cal D}^{'},\\
&\rho_{n}^{t}\n v_{n}\rightarrow \rho^{t}\n v\;\;\;\mbox{in}\;\;{\cal D}^{'},\\
\end{aligned}
$$
\end{lemme}
{\bf Proof:} Let $\phi$ be a test function, then:
$$
\begin{aligned}
&\int \rho_{n}\n v_{n}\phi dxdt=-\int \rho_{n}v_{n}\n\phi dxdt+\int v_{n}\n \rho_{n}\,\phi dxdt\\
&\hspace{2cm}=-\int\sqrt{\rho_{n}}(\sqrt{\rho_{n}}v_{n})\n\phi dxdt+\int \sqrt{\rho_{n}}v_{n}\n(2\sqrt{\rho_{n}})\phi dxdt.\\
\end{aligned}
$$
Thanks to lemma \ref{lemme1}, we know that $\sqrt{\rho_{n}}$ is bounded in $L^{\infty}(0,T;L^{6}_{loc}(\Omega))$ and that $\sqrt{\rho_{n}}$ converges almost everywhere to $\sqrt{\rho}$ (defined to be zero on the vacuum set). Therefore it converges strongly in $L^{2}_{loc}((0,T)\times\Omega)$. This point is enough to prove the convergence of the first term.\\
Next as $\sqrt{\rho_{n}}$ converges strongly to $\sqrt{\rho}$ in $L^{2}_{loc}((0,T)\times\Om)$, it follows that:
$$\n\sqrt{\rho_{n}}\rightarrow\n\sqrt{\rho}\;\;\;L^{2}_{loc}((0,T)\times\Omega)-\mbox{weak}.$$
And we conclude for the second term as $\sqrt{\rho_{n}}v_{n}$ converges strongly in $L^{2}_{loc}((0,T)\times\Omega)$.
\subsubsection*{Convergence of $\sqrt{\rho_{n}}v_{n}$ to $\sqrt{\rho}u+\frac{2\kappa}{\mu}\n\sqrt{\rho}$}
\begin{lemme}
We have the following property.
$$
\sqrt{\rho_{n}}v_{n}\rightarrow \sqrt{\rho}u+\frac{2\kappa}{\mu}\n\sqrt{\rho}\;\;\mbox{a.e}.
$$
\end{lemme}
{\bf Proof:} As:
$$\D\rho_{n}=\sqrt{\rho_{n}}(\sqrt{\rho_{n}}\D\ln\rho_{n})+4|\n\sqrt{\rho_{n}}|^{2},$$
we can show by (\ref{27}) that $\D\rho_{n}$ is uniformly bounded in $L_{T}^{\infty}(L^{1}(\Om))+L^{2}_{T}(L^{\frac{3}{2}}(\Om))$. 
By Aubin Lions theorem, we deduce that $\n\rho_{n}$ converges a.e to $\n\rho$. It concludes the proof of the lemma.
\subsection{Stability for theorem \ref{2.2}}
The proof follows the same line than for theorem \ref{2.1} and than \cite{fMV}.
\section{Proof of theorem \ref{theo4}}
\label{section6}
For the result of existence of strong solution in the theorem \ref{theo4}, we refer to \cite{Hprepa}. We now want to concentrate on the result of blow-up of the theorem \ref{ftheo4}. To do this, we shall point out new regularizing effects and gain of integrability on the density and the effective velocity $v$. Indeed the key of the proof is to obtain enough integrability on $v$ to derive regularizing effects on the density via the first equation of system \ref{3systeme2} (which is a heat equation with a remainder term which depends on $v$, i.e ${\rm div}(\rho v)$).
\subsubsection*{Gain of integrability on $v$}
As in \cite{regular}, we now want to obtain additional information on the integrability of $v$, and more precisely we would like to show that $\rho^{\frac{1}{p}}v$ is in any $L^{\infty}_{T}(L^{p}(\T^{N})$ with $1<p<+\infty$.  To do it, we multiply the momentum equation of (\ref{3systeme2}) by $v|v|^{p-2}$ and integrate over $\R^{N}$, we obtain then:
\begin{equation}
\begin{aligned}
&\frac{1}{p}\int_{\T^{N}}\rho\p_{t}(|v|^{p})dx+\int_{\T^{N}}\rho u\cdot\n(\frac{|v|^{p}}{p})dx
+\int_{\T^{N}} \rho|v|^{p-2}|\n v|^{2}dx\\
&\hspace{1,5cm}+(p-2)\int_{\T^{N}} \rho\sum_{i,j,k}v_{j}v_{k}\p_{i}v_{j}\p_{i}v_{k}|v|^{p-4}dx+\int_{\T^{N}} |v|^{p-2}v\cdot\n a\rho dx=0.
\end{aligned}
\label{in1A}
\end{equation}
Next we observe that:
$$\sum_{i,j,k}v_{j}v_{k}\p_{i}v_{j}\p_{i}v_{k}=\sum_{i}(\sum_{j}v_{j}\p_{i}v_{j})^{2}=\sum_{i}\frac{1}{2}\p_{i}(|v|^{2}).$$
We get then as ${\rm div}(\rho u)=-\p_{t}\rho$ and by using (\ref{in1A}):
\begin{equation}
\begin{aligned}
&\frac{1}{p}\int_{\T^{N}}\p_{t}(\rho|v|^{p})dx+\int_{\T^{N}} \rho|v|^{p-2}|\n v|^{2}dx+(p-2)\int_{\T^{N}} \rho(\sum_{i}\p_{i}(|v|^{2})^{2}|v|^{p-4}dx\\
&\hspace{8cm}+\int_{\T^{N}} |v|^{p-2}v\cdot\n a\rho dx=0.
\end{aligned}
\label{in1A1}
\end{equation}
We have then by integrating over $(0,t)$ with $0<t\leq T$:
\begin{equation}
\begin{aligned}
&\frac{1}{p}\int_{\T^{N}}(\rho|v|^{p})(t,x)dx+\int^{t}_{0}\int_{\T^{N}} \rho|v|^{p-2}|\n v|^{2}(t,x)dtdx\\
&+(p-2)\int^{t}_{0} \int_{\T^{N}}\rho(\sum_{i}\p_{i}(|v|^{2})^{2}|v|^{p-4}(t,x)dtdx\leq \frac{1}{p}\int_{\T^{N}}(\rho_{0}|v_{0}|^{p})(x)dx\\
&\hspace{7cm}+|\int^{t}_{0} \int_{\T^{N}}|v|^{p-2}v\cdot\n a\rho(t,x)dtdx|.
\end{aligned}
\label{in1A1}
\end{equation}
We now want to take the $\sup$ of the previous estimate on $(0,T)$, we have then:
\begin{equation}
\begin{aligned}
&\sup_{t\in(0,T)}\big(\frac{1}{p}\int_{\T^{N}}(\rho|v|^{p})(t,x)dx+\int^{t}_{0}\int_{\T^{N}} \rho|v|^{p-2}|\n v|^{2}(s,x)dsdx\\
&+(p-2)\int^{t}_{0} \int_{\T^{N}}\rho(\sum_{i}\p_{i}(|v|^{2})^{2}|v|^{p-4}(s,x)dtdx\big)\leq \frac{1}{p}\int_{\T^{N}}(\rho_{0}|v_{0}|^{p})(x)dx\\
&\hspace{5,5cm}+\sup_{t\in(0,T)}|\int^{t}_{0} \int_{\T^{N}}|v|^{p-2}v\cdot\n a\rho(s,x)dtdx|.
\end{aligned}
\label{in1A1a}
\end{equation}
By integration by parts we have:
$$\int^{t}_{0} \int_{\T^{N}}|v|^{p-2}v\cdot\n a\rho(t,x)dtdx=-\int^{t}_{0} \int_{\T^{N}}{\rm div}(|v|^{p-2}v)a\rho(t,x)dtdx.$$
We have then:
$${\rm div}(|v|^{p-2}v)=|v|^{p-2}{\rm div}(v)+(p-2)|v|^{p-4}v\cdot(v\cdot\n v).$$
We now want to bound $\sup_{t\in(0,T)}|\int^{t}_{0} \int_{\T^{N}}|v|^{p-2}v\cdot\n a\rho(t,x)dtdx|$ in (\ref{in1A1}). To do this, we would like to apply H\"older's inequalities, indeed we have:
$$|\rho|v|^{p-4}v\cdot(v\cdot\n v)|\leq C|\sqrt{\rho}|v|^{\frac{p}{2}-1}\n v|\times|\rho^{\frac{1}{2}-\frac{1}{p}}|v|^{\frac{p}{2}-1}|\times|\rho^{\frac{1}{p}}|.$$
By H\"older«s inequalities we can prove that $|\rho|v|^{p-4}v\cdot(v\cdot\n v)|$ is in $L^{1}_{t}(L^{1}(\T^{N})$ with:
\begin{equation}
\begin{aligned}
&\sup_{t\in(0,T)}\int^{t}_{0}\int_{\T^{N}}|\rho|v|^{p-4}v\cdot(v\cdot\n v)|dxds\leq \sup_{t\in(0,T)} \big(\sqrt{t}\|\sqrt{\rho}|v|^{\frac{p}{2}-1}\n v\|_{L^{2}_{t}(L^{2}(\T^{N}))}\\
&\hspace{5cm}\times\|\rho^{\frac{1}{2}-\frac{1}{p}}|v|^{\frac{p}{2}-1}|\|_{L^{\infty}_{t}(L^{\frac{2p}{p-2}}(\T^{N}))}\|\rho^{\frac{1}{p}}\|_{L^{\infty}_{t}(L^{p}(\T^{N}))}\big),\\[2mm]
&\leq \sup_{t\in(0,T)}\big(\sqrt{t}\|\sqrt{\rho}|v|^{\frac{p}{2}-1}\n v\|_{L^{2}_{t}(L^{2}(\T^{N}))}\sup_{s\in(0,t)}\big(\int_{\T^{N}}\rho|v||^{p}(x)dx\big)^{\frac{1}{2}-\frac{1}{p}}\|\rho^{\frac{1}{p}}\|_{L^{\infty}_{t}(L^{p}(\T^{N}))}\big).
\end{aligned}
\label{1int2}
\end{equation}
Now we set:
$$
\begin{aligned}
&A(t)=\big(\frac{1}{p}\int_{\T^{N}}(\rho|v|^{p})(t,x)dx+\int^{t}_{0}\int_{\T^{N}} \rho|v|^{p-2}|\n v|^{2}(s,x)dsdx\\
&\hspace{6cm}+(p-2)\int^{t}_{0} \int_{\T^{N}}\rho(\sum_{i}\p_{i}(|v|^{2})^{2}|v|^{p-4}(s,x)dtdx\big).
\end{aligned}
$$
We have then from (\ref{in1A1a}) and (\ref{1int2}):
$$\sup_{t\in(0,T)}A(t)\leq C \big(1+\sqrt{T}(\sup_{t\in(0,T)}A(t))^{1-\frac{1}{p}}\|\rho\|_{L^{\infty}_{T}(L^{1}(\T^{N}))}^{\frac{1}{p}}\big)
$$
We deduce then that $\sup_{t\in(0,T)}A(t)$ is finite. In particular, we obtain that $\rho^{\frac{1}{p}}v$ is in any $L^{\infty}_{T}(L^{p}(\T^{N}))$ for $1<p<+\infty$.
\subsubsection*{Blow-up criterion}
We have then obtain that  for any $1\leq p<+\infty$, $\rho^{\frac{1}{p}}v$ belongs to $L^{\infty}(L^{p})$. We now want to take into account this information to obtain regularizing effects on the density via the first equation of (\ref{3systeme2}) (which is a heat equation):
\begin{equation}
\p_{t}\rho-\frac{\kappa}{\mu}\D \rho=-{\rm div}(\rho v).
\label{charegul}
\end{equation}
Our goal is to transfer the information on the integrability of $v$ (which is a \textit{subscaling} estimate) on the density $\rho$. More precisely we have by proposition \ref{5chaleur} for any $1\leq p<+\infty$:
\begin{equation}
\|\rho\|_{\widetilde{L}_{T}^{\infty}(B^{1}_{p,\infty})}\leq C(\|\rho_{0}\|_{\widetilde{L}^{\infty}(B^{1}_{p,\infty})}+\|\rho v\|_{\widetilde{L}_{T}^{\infty}(B^{0}_{p,\infty})}).
\label{incha}
\end{equation}
 We now need to prove that $\rho v$ is in $\widetilde{L}_{T}^{\infty}(B^{0}_{p,\infty})$. As $\rho^{\frac{1}{p}}v$ is in $L^{\infty}_{T}(L^{p})$, we have then:
\begin{equation}
\begin{aligned}
\|(\rho^{1-\frac{1}{p}})\rho^{\frac{1}{p}}v\|_{L_{T}^{\infty}(L^{p}(\T^{N}))}&\leq\|\rho^{\frac{1}{p}}v\|_{L_{T}^{\infty}(L^{p})}\|\rho^{1-\frac{1}{p}}\|_{L^{\infty}_{T}(L^{\infty}(\T^{N}))},\\
&\leq  C\|\rho^{\frac{1}{p}}v\|_{L_{T}^{\infty}(L^{p})}\|\rho\|^{1-\frac{1}{p}}_{L^{\infty}_{T}(L^{\infty}(\T^{N}))}.
\end{aligned}
\label{prod1}
\end{equation}
By injecting (\ref{prod1})  in (\ref{incha}) and using the fact that $L^{p}(\T^{N})$ is embedded in $B^{0}_{p,\infty}$, we obtain that:
\begin{equation}
\|\rho\|_{\widetilde{L}_{T}^{\infty}(\dot{B}^{1}_{p,\infty})}\leq C(\|\rho_{0}\|_{\widetilde{L}^{\infty}(\dot{B}^{1}_{p,\infty})}+\|\rho^{\frac{1}{p}} v\|_{L_{T}^{\infty}(L^{p}(\T^{N}))}\|\rho\|^{1-\frac{1}{p}}_{L^{\infty}_{T}(L^{\infty}(\T^{N}))} \big).
\label{incha12}
\end{equation}
By Besov embedding we know that $\dot{B}^{1}_{p,\infty})$ is embedded in $L^{\infty}(\T^{N})$ for $p$ large enough. We obtain by (\ref{incha12}) that:
\begin{equation}
\|\rho\|_{\widetilde{L}_{T}^{\infty}(\dot{B}^{1}_{p,\infty})}\leq C(\|\rho_{0}\|_{\widetilde{L}^{\infty}(B^{1}_{p,\infty})}+\|\rho^{\frac{1}{p}} v\|_{L_{T}^{\infty}(L^{p}(\T^{N}))}\|\rho\|^{1-\frac{1}{p}}_{L^{\infty}_{T}(B^{1}_{p,\infty}))} \big).
\label{incha12}
\end{equation}
We conclude then by bootstrap that $\rho$ is in $L_{T}^{\infty}(B^{1}_{p,\infty})$ for $p$ large enough (in fact $p>N$).
\subsubsection*{Estimates on $\frac{1}{\rho}$ if we control $\frac{1}{\rho^{\e}}$ in  $L^{\infty}_{T}(L^{1}(\T^{N}))$}
We now would like to obtain estimates on $\frac{1}{\rho}$, to do this we will consider the first equation of (\ref{3systeme2}) that we multiply by $-\frac{1}{\rho^{p}}$ with $p\geq 2$:
\begin{equation}
\frac{1}{p-1}\p_{t}(\frac{1}{\rho^{p-1}})-\frac{\kappa}{\mu(p-1)}\D(\frac{1}{\rho^{p-1}})+\frac{4p\kappa}{\mu(p-1)^{2}}|\n\big(\frac{1}{\rho^{\frac{p}{2}-\frac{1}{2}}}\big)|^{2}=-\frac{1}{\rho^{p}}{\rm div}(\rho v),
\label{vide1}
\end{equation}
because we have:
$$\frac{\kappa}{\mu\rho^{p}}\D\rho=-\frac{\kappa}{\mu(p-1)}\D(\frac{1}{\rho^{p-1}})+\frac{4p\kappa}{\mu(p-1)^{2}}|\n\big(\frac{1}{\rho^{\frac{p}{2}-\frac{1}{2}}}\big)|^{2}.$$
Our goal is now to integrate over $(0,t)\times\T^{N}$ the equality (\ref{vide1}), we obtain then:
\begin{equation}
\begin{aligned}
&\frac{1}{p-1}\int_{\T^{N}}\frac{1}{\rho^{p-1}(t,x)}dx+\frac{4p\kappa}{\mu(p-1)^{2}}\int^{t}_{0}\int_{\T^{N}}|\n\big(\frac{1}{\rho^{\frac{p}{2}-\frac{1}{2}}(t,x)}\big)|^{2}dtdx\\
&\hspace{3cm}\leq \frac{1}{p-1}\int_{\T^{N}}\frac{1}{\rho_{0}^{p-1}(x)}dx+\big|\int^{t}_{0}\int_{\T^{N}}
\frac{1}{\rho^{p}}{\rm div}(\rho v)(t,x)dtdx\big|.
\end{aligned}
\label{vide2}
\end{equation}
We now pass to the $\sup$ on $(0,T)$, and we have:
\begin{equation}
\begin{aligned}
&\sup_{t\in(0,T)}\big(\frac{1}{p-1}\int_{\T^{N}}\frac{1}{\rho^{p-1}(t,x)}dx+
\frac{4p\kappa}{\mu(p-1)^{2}}\int^{t}_{0}\int_{\T^{N}}|\n\big(\frac{1}{\rho^{\frac{p}{2}-\frac{1}{2}}(t,x)}\big)|^{2}dtdx\big)\\
&\hspace{2cm}\leq \frac{1}{p-1}\int_{\T^{N}}\frac{1}{\rho_{0}^{p-1}(x)}dx+\sup_{t\in(0,T)}\big( \big|\int^{t}_{0}\int_{\T^{N}}
\frac{1}{\rho^{p}}{\rm div}(\rho v)(t,x)dtdx\big|\big).
\end{aligned}
\label{vide3}
\end{equation}
By integration by parts we get:
$$\big|\int^{t}_{0}\int_{\T^{N}}
\frac{1}{\rho^{p}}{\rm div}(\rho v)(t,x)dtdx\big|=|p\sum_{i}\int^{t}_{0}\int_{\T^{N}}
\frac{1}{\rho^{p}}\p_{i}\rho\,v^{i}(t,x)dtdx|.$$
We recall that:
\begin{equation}
|\frac{1}{\rho^{p}}\p_{i}\rho\,v^{i} |\big(t,x\big) =|\frac{1}{\rho^{\frac{p}{2}+\frac{1}{2}}}\p_{i}\rho|\times|\frac{1}{\rho^{\frac{p}{2}-\frac{1}{2}}}|\times|v^{i}|\big(t,x\big).
\label{aHolder}
\end{equation}
In the sequel we only shall deal with the case $N=3$, the other cases are similar to treat. Now as on the left side of (\ref{vide3}), we have information on $\frac{1}{\rho^{\frac{p}{2}-\frac{1}{2}}}$ in $L^{\infty}_{T}(L^{2}(\T^{N}))$ and on $\n(\frac{1}{\rho^{\frac{p}{2}-\frac{1}{2}}})$ in $L^{2}_{T}(L^{2}(\T^{N}))$ we can assume by Sobolev embedding and interpolation that:
\begin{equation}
\frac{1}{\rho^{\frac{p}{2}-\frac{1}{2}}}\in L^{k}_{T}(L^{q}(\T^{N}))\;\;\;\mbox{with}\frac{1}{q}+\frac{2}{3k}=\frac{1}{2}.
\label{aainterpo1}
\end{equation}
We have then:
\begin{equation}
\begin{aligned}
\|\frac{1}{\rho^{\frac{p}{2}-\frac{1}{2}}}\|_{L^{k}_{T}(L^{q}(\T^{N}))}&\leq
\|\frac{1}{\rho^{\frac{p}{2}-\frac{1}{2}}}\|_{L^{\infty}_{T}(L^{2}(\T^{N}))}^{1-\frac{2}{p}}
\|\frac{1}{\rho^{\frac{p}{2}-\frac{1}{2}}}\|_{L^{2}_{T}(L^{6}(\T^{N}))}^{\frac{2}{p}},\\
&\leq \|\frac{1}{\rho^{\frac{p}{2}-\frac{1}{2}}}\|_{L^{\infty}_{T}(L^{2}(\T^{N}))}^{1-\frac{2}{p}}\big(
\|\n\big(\frac{1}{\rho^{\frac{p}{2}-\frac{1}{2}}}\big\|_{L^{2}_{T}(L^{2}(\T^{N}))}+\sqrt{T}\|\frac{1}{\rho^{\frac{p}{2}-\frac{1}{2}}}\|
_{L^{\infty}_{T}(L^{2}(\T^{N}))}
\big)^{\frac{2}{p}}.
\end{aligned}
\label{bHolder}
\end{equation}
We now want to take advantage of the fact that $\frac{1}{\rho^{\e}}$ belongs to $L^{\infty}_{T}(L^{1}(\T^{N}))$ for $\e$ arbitrary small. Indeed the idea is now to use interpolqtion results in order to absorb the term on the right side of (\ref{vide3}) by the left side. More precisely we have by (\ref{aainterpo1}) that:
\begin{equation}
\frac{1}{\rho^{\e}}\in L^{\infty}_{T}(L^{1}(\T^{N}))\;\;\;\mbox{and}\;\;\;\frac{1}{\rho^{\e}}\in L^{k(\frac{p-1}{2\e})}_{T}(L^{q(\frac{p-1}{2\e})}(\T^{N})),
\label{aainterpo12}
\end{equation}
with:
\begin{equation}
\|\frac{1}{\rho^{\e}}\|_{L^{k(\frac{p-1}{2\e})}_{T}(L^{q(\frac{p-1}{2\e})}(\T^{N}))}=\|\frac{1}{\rho^{\frac{p}{2}-\frac{1}{2}}}\|_{L^{k}_{T}(L^{q}(\T^{N}))}
^{\frac{2\e}{(p-1)}},
\label{aainterpo13}
\end{equation}
Now by interpolation we can show that:
\begin{equation}
\|\frac{1}{\rho^{\e}}\|_{L^{\alpha}_{T}(L^{\beta}(\T^{N}))}\leq \|\frac{1}{\rho^{\e}}\|^{1-\theta}_{L^{\infty}_{T}(L^{1}(\T^{N}))} \|\frac{1}{\rho^{\e}}\|^{\theta}_{L^{k(\frac{p-1}{2\e})}_{T}(L^{q(\frac{p-1}{2\e})}(\T^{N}))},
\label{aainterpo14}
\end{equation}
with $\frac{1}{\alpha}=\frac{2\e\theta}{k(p-1)}$ and $\frac{1}{\beta}=(1-\theta)+\frac{2\e\theta}{q(p-1)}$ and $0\leq\theta\leq 1$.\\
We check easily that:
\begin{equation}
\|\frac{1}{\rho^{\e}}\|_{L^{\alpha}_{T}(L^{\beta}(\T^{N}))}=\|\frac{1}{\rho^{\frac{p}{2}-\frac{1}{2}}}\|
^{\frac{p-1}{2\e}}_{L^{\frac{2\e\alpha}{p-1}}_{T}(L^{\frac{2\e\beta}{p-1}}(\T^{N}))}.
\label{aainterpo15}
\end{equation}
From (\ref{aainterpo13}), (\ref{aainterpo14}) and (\ref{aainterpo15}), we obtain finally that:
\begin{equation}
\|\frac{1}{\rho^{\frac{p}{2}-\frac{1}{2}}}\|
^{\frac{p-1}{2\e}}_{L^{\frac{2\e\alpha}{p-1}}_{T}(L^{\frac{2\e\beta}{p-1}}(\T^{N}))}\leq \|\frac{1}{\rho^{\e}}\|^{1-\theta}_{L^{\infty}_{T}(L^{1}(\T^{N}))}
\|\frac{1}{\rho^{\frac{p}{2}-\frac{1}{2}}}\|_{L^{k}_{T}(L^{q}(\T^{N}))},
^{\frac{2\e\theta}{(p-1)}}
\label{aainterpo16}
\end{equation}
or:
\begin{equation}
\|\frac{1}{\rho^{\frac{p}{2}-\frac{1}{2}}}\|_{L^{\frac{2\e\alpha}{p-1}}_{T}(L^{\frac{2\e\beta}{p-1}}(\T^{N}))}\leq \|\frac{1}{\rho^{\e}}\|^{\frac{2\e}{(p-1)}(1-\theta)}_{L^{\infty}_{T}(L^{1}(\T^{N}))}
\|\frac{1}{\rho^{\frac{p}{2}-\frac{1}{2}}}\|_{L^{k}_{T}(L^{q}(\T^{N}))}
^{(\frac{2\e}{(p-1)})^{2}\theta}.
\label{aainterpo17}
\end{equation}
We recall that:
\begin{equation}
\frac{2\e\alpha}{p-1}=\frac{k}{\theta}\;\;\;\mbox{and}\;\;\;\frac{2\e\beta}{p-1}=\frac{2\e q}{(1-\theta)(p-1)q+2\e\theta}.
\label{acoeffs}
\end{equation}
By using (\ref{aHolder}), (\ref{bHolder}) and (\ref{aainterpo17}) we obtain by H\"older'sinequelities:
\begin{equation}
\begin{aligned}
&\big|\int^{t}_{0}\int_{\T^{N}}
\frac{1}{\rho^{p}}{\rm div}(\rho v)(t,x)dtdx\big|\leq C\|\frac{1}{\rho^{\frac{p}{2}+\frac{1}{2}}}\p_{i}\rho\|_{L^{2}_{T}(L^{2}(\T^{N}))}
\\
&\hspace{3cm}\times\|\frac{1}{\rho^{\frac{p}{2}-\frac{1}{2}}}\|_{L^{k(1+\alpha_{1})}_{T}(L^{q(1+\alpha_{2})}(\T^{N}))}\|v\|_{L^{k^{1}}_{T}(L^{q^{1}}(\T^{N}))}.
\end{aligned}
\label{CHolder}
\end{equation}
with $\alpha_{1}>0$, $\alpha_{2}>0$ arbitrary small, $\theta$ depending on  $\alpha_{1}>0$ and $\alpha_{2}$ (such that $0<\theta<1$ and $\theta$ goes to $1$ when $\alpha_{1}$ and $\alpha_{2}$ go to $0$) and $k^{1}$, $k^{2}$ such that:
$$\frac{1}{2}+\frac{1}{k(1+\alpha_{1})}+\frac{1}{k^{1}}=1\;\;\;\mbox{and}\;\;\;\frac{1}{2}+\frac{1}{q(1+\alpha_{2})}+\frac{1}{k^{2}}=1.$$
From (\ref{bHolder}), (\ref{CHolder}) and (\ref{aainterpo17}) we have finaly:
\begin{equation}
\begin{aligned}
&\big|\int^{t}_{0}\int_{\T^{N}}
\frac{1}{\rho^{p}}{\rm div}(\rho v)(t,x)dtdx\big|\leq C\big(\sup_{t\in(0,T)}B(t)\big)^{1+(\frac{2\e}{(p-1)})^{2}\theta}\\
&\hspace{6cm}\|\frac{1}{\rho^{\e}}\|^{\frac{2\e}{(p-1)}(1-\theta)}_{L^{\infty}_{T}(L^{1}(\T^{N}))}\|v\|_{L^{k^{1}}_{T}(L^{q^{1}}(\T^{N}))},
\end{aligned}
\end{equation}
with:
$$B(t)=\big(\frac{1}{p-1}\int_{\T^{N}}\frac{1}{\rho^{p-1}(t,x)}dx+
\frac{4p\kappa}{\mu(p-1)^{2}}\int^{t}_{0}\int_{\T^{N}}|\n\big(\frac{1}{\rho^{\frac{p}{2}-\frac{1}{2}}(t,x)}\big)|^{2}dtdx\big).$$
When $(p-1)>2\e$ we can absorb the term $|\int^{t}_{0}\int_{\T^{N}}
\frac{1}{\rho^{p}}{\rm div}(\rho v)(t,x)dtdx\big|$ because from the previous inequality, we have:
$$\big(\sup_{t\in(0,T)}B(t)\big)^{2}\leq C\big(1+\big(\sup_{t\in(0,T)}B(t)\big)^{1+(\frac{2\e}{(p-1)})^{2}\theta})\|v\|_{L^{k^{1}}_{T}(L^{q^{1}}(\T^{N}))}\big),$$
with $1+(\frac{2\e}{(p-1)})^{2}\theta<2$. Indeed we can show easily that $v$ is bounded in $L^{k^{1}}_{T}(L^{q^{1}}(\T^{N}))$.
To do this, we recall that $\rho^{\frac{1}{p}}v$ is bounded in $L^{\infty}_{T}(L^{p}(\T^{N}))$ for $p$ enough big. As we control $\frac{1}{\rho^{\e}}$ in 
$L^{\infty}_{T}(L^{1}(\T^{N}))$ and that:
$$v=\frac{1}{\rho^{\frac{1}{p}}}\rho^{\frac{1}{p}}v,$$
$V$ is then bounded in $L^{\infty}_{T}(L^{\frac{p}{2}}(\T^{N}))$ for $p$ enough big and then $v$ is bounded in $L^{k^{1}}_{T}(L^{q^{1}}(\T^{N}))$.\\
We have then proved that for any $p>0$, $\frac{1}{\rho^{p}}$ is bounded in $L^{\infty}_{T}(L^{1}(\T^{N}))$.
\subsubsection*{Estimates on $\frac{1}{\rho}$ if we control $v$ in  $L^{p}_{T}(L^{q}(\T^{N}))$ with $\frac{1}{p}+\frac{N}{2q}=\frac{1}{2}$ and $1\leq p<+\infty$}
In this case, it suffices only to bound $\big|\int^{t}_{0}\int_{\T^{N}}
\frac{1}{\rho^{p}}{\rm div}(\rho v)(t,x)dtdx\big|$ as follows:
$$
\begin{aligned}
\big|\int^{t}_{0}\int_{\T^{N}}
\frac{1}{\rho^{p}}{\rm div}(\rho v)(t,x)dtdx\big|\leq C \|\frac{1}{\rho^{\frac{p}{2}+\frac{1}{2}}}\p_{i}\rho\|_{L^{2}_{T}(L^{2}(\T^{N}))}\|\frac{1}{\rho^{\frac{p}{2}-\frac{1}{2}}}\|_{L^{p^{'}}_{T}(L^{q^{'}}(\T^{N}))}
\|v\|_{L^{p}_{T}(L^{q}(\T^{N}))}.
\end{aligned}
$$
And by Sobolev embedding we have:
$$
\begin{aligned}
\big|\int^{t}_{0}\int_{\T^{N}}
\frac{1}{\rho^{p}}{\rm div}(\rho v)(t,x)dtdx\big|\leq C \big(\sup_{t\in(0,T)}B(t)\big)^{2}
\|v\|_{L^{p}_{T}(L^{q}(\T^{N}))}.
\end{aligned}
$$
This term can be absorbed in (\ref{vide3}) because if $T$ is small enough, then $\|v\|_{L^{p}_{T}(L^{q}(\T^{N}))}$ is small enough to do a bootstrap. To obtain the result on a general interval $(0,T)$ it suffices to decompose the interval in a sum of small intervals and to apply the previous idea by iteration.
\subsubsection*{Control of $\n\ln\rho$ in $L^{\infty}(B^{0}_{N+\e,\infty})$}
To do this it suffice to recall that $\n\sqrt{\rho}$ is in $L^{\infty}_{T}(L^{2}(\T^{N}))$, and as $\frac{1}{\rho^{\alpha}}$ is in $L^{\infty}_{T}(L^{1}(\T^{N}))$ for all $0<\alpha$ then by Sobolev embedding we can show easily that $\frac{1}{\rho}$ is in $L^{\infty}_{T}(B^{0}_{N+\e,1})$. Since we have seen that $\n\rho$ is in $L^{\infty}_{T}(B^{0}_{p,\infty})$ for any $p>1$, by paraproduct we obtain that:
$$\n\ln\rho\in L^{\infty}(B^{0}_{N+\e,\infty}).$$
Furthermore we have obtained that $\rho^{\frac{1}{p}}v$ belongs to $L^{\infty}(L^{p})$ for any $1\leq p<+\infty$ and that $\n\ln\rho\in L^{\infty}_{T}(B^{0}_{N+\e,\infty})$. As $u=v-\frac{\kappa}{\mu}\n\ln\rho$, we obtain easily by following the same lines that $u$ is in $L^{\infty}_{T}(B^{0}_{N+\e,\infty}$. To summarize what we have obtained, we have:
\begin{equation}
u\in L^{\infty}_{T}(B^{0}_{N+\e,\infty})\;\;\;\mbox{and}\;\;\;\ln\rho\in L^{\infty}_{T}(\dot{B}^{1}_{N+\e,\infty}).
\label{hyperimp}
\end{equation}
Easily we can prove that if $u_{0}\in B^{0}_{N+\e,\infty})$ and $\ln\rho_{0}\in B^{1}_{N+\e,\infty}$ then the system (\ref{3systeme2}) has a strong solution $(\rho,u)$ on $(0,T^{'})$ with:
$$T^{'}\geq \frac{C}{(1+\|u_{0}\|_{B^{0}_{N+\e,\infty}}+\|\ln\rho_{0}\|_{B^{1}_{N+\e,\infty}})^{\beta}}.$$
We refer to \cite{Hprepa} for the proof. This is an easy consequence of the fact that the initial data are choose subcritical.\\
It means that there exists a time $T^{'}\geq c>0$, where $c$ depends only on the physical coefficients and of subcritical initial data. We can construct by theorem \ref{ftheo1} a solution $(\rho_{1},u_{1})$ on $(T-\alpha,T-\alpha+T^{'})$ with initial data $(\rho(T-\alpha), u(T-\alpha))$ (here $\alpha<T^{'}$). The only difficulty is to prove that on $(T-\alpha,T)$ we have:
$$(\rho_{1},u_{1})=(\rho,u).$$
To do this, it suffices only to use the uniqueness part of theorem 1 in \cite{Hprepa}. It concludes the proof of theorem \ref{theo4}.
\hfill {$\Box$}
\section{Appendix}
In this appendix, we only want to detail the computation on the Korteweg tensor.
\begin{lemme}
$${\rm div}K=\kappa{\rm div}(\rho\n\n\ln\rho)=\kappa{\rm div}(\rho D(\n\ln\rho)).$$
\end{lemme}
{\bf Proof:} By calculus, we obtain then:
\begin{equation}
\begin{aligned}
({\rm div}K)_{j}&=\big(\n\D\rho-{\rm div}(\frac{1}{\rho}\n\rho\otimes\n\rho)\big)_{j},\\
&=\p_{j}\D\rho-\frac{1}{\rho}\D\rho\,\p_{j}\rho-\frac{1}{2\rho}\p_{j}|\n\rho|^{2}+\frac{1}{\rho^{2}}|\n\rho|^{2}\p_{j}\rho,\\
\end{aligned}
\label{princip}
\end{equation}
Next we have:
$$\D\rho=\rho\D\ln\rho+\frac{1}{\rho}|\n\rho|^{2}.$$
We have then:
\begin{equation}
\begin{aligned}
\p_{j}\D\rho-\frac{1}{\rho}\D\rho\p_{j}\rho&=\p_{j}(\rho\,\D\ln\rho+\frac{1}{\rho}|\n\rho|^{2})
-\D\ln\rho\,\p_{j}\rho-\frac{1}{\rho^{2}}|\n\rho|^{2}\p_{j}\rho,\\
&=\rho\p_{j}\D\ln\rho+\frac{1}{\rho}\p_{j}(|\n\rho|^{2})-\frac{2}{\rho^{2}}|\n\rho|^{2}\p_{j}\rho.
\end{aligned}
\label{prin1}
\end{equation}
Putting the expression of (\ref{prin1}) in (\ref{princip}), we obtain:
\begin{equation}
\begin{aligned}
({\rm div}K)_{j}&=\p_{j}\D\rho+\frac{1}{2\rho}\p_{j}(|\n\rho|^{2})-\frac{1}{\rho^{2}}|\n\rho|^{2}\p_{j}\rho.
\end{aligned}
\label{principa}
\end{equation}
Next by calculus, we have:
\begin{equation}
\begin{aligned}
\frac{1}{2\rho}\p_{j}(|\n\rho|^{2})-\frac{1}{\rho^{2}}|\n\rho|^{2}\p_{j}\rho&=\sum_{i}(\p_{i}\ln\rho\p_{ij}\rho-(\p_{i}\ln\rho)^{2}\p_{j}\rho),\\
&=\sum_{i}\p_{i}\ln\rho\,\rho\p_{i,j}\ln\rho,\\
&=\frac{\rho}{2}\n(|\ln\rho|^{2})_{j}.
\end{aligned}
\label{imp}
\end{equation}
Finally by using (\ref{imp}) and (\ref{principa}), we obtain:
$$
\begin{aligned}
{\rm div}K=\rho(\n\D(\ln\rho)+\frac{\rho}{2}\n(|\n\ln\rho|^{2})).
\label{cap}
\end{aligned}
$$
We now want to prove that we can rewrite (\ref{cap}) under the form of a viscosity tensor. To see this, we have:
$$
\begin{aligned}
{\rm div}(\rho\n(\n \ln\rho))_{j}&=\sum_{i}\p_{i}(\rho\p_{ij}\ln\rho),\\
&=\sum_{i}[\p_{i}\rho\p_{ij}\ln\rho+\rho\p_{iij}\ln\rho],\\
&=\rho(\D\n \ln\rho)_{j}+\sum_{i}\rho\p_{i}\ln\rho\p_{j}\p_{i}\ln\rho),\\
&=\rho(\D\n \ln\rho)_{j}+\frac{\rho}{2}(\n(|\n\ln\rho|^{2}))_{j},\\
&={\rm div}K.
\end{aligned}
$$
We have then:
$${\rm div}K=\kappa{\rm div}(\rho\n\n\ln\rho)=\kappa{\rm div}(\rho D(\n\ln\rho)).$$

\end{document}